\documentclass[10pt]{amsart}
\usepackage{amsmath,amsthm,amssymb}

\usepackage{amsmath,amsthm,graphicx}
\usepackage[colorlinks=true,linkcolor=blue,urlcolor=blue,citecolor=blue]{hyperref}
\usepackage[all]{xy}

  \def\doi#1{%
    \href{https://doi.org/#1}{doi:#1}}%

\newtheorem{teo}{Theorem}[section]
\newtheorem{rmk}{Remark}[section]

\newtheorem{prp}{Proposition}[section]
\newtheorem{cor}{Corollary}[section]

\usepackage{mathtools}
\DeclarePairedDelimiter\floor{\lfloor}{\rfloor}
\newcommand{\ep}{\varepsilon_{\mathcal{T}}}
\newcommand{\R}{\mathbb R}

\newcommand{\U}{\,\mathcal U\,}

\DeclareMathOperator*{\argmax}{arg\,max}
\DeclareMathOperator*{\argmin}{arg\,min}

\def\cGo{\overline{\Omega}}     
\def\fGo{\partial\Omega }         
\def\Ux{\U_x}

\begin{document}

\newtheorem{theorem}{Theorem}[section]
\newtheorem{lemma}[theorem]{Lemma}
\newtheorem{proposition}[theorem]{Proposition}
\newtheorem{corollary}[theorem]{Corollary}

\theoremstyle{definition}
\newtheorem{definition}[theorem]{Definition}
\newtheorem{remark}[theorem]{Remark}
\newtheorem{example}[theorem]{Example}
\newtheorem{problem}[theorem]{Problem}

\numberwithin{equation}{section}

\renewcommand{\subjclassname}{%
\textup{2010} Mathematics Subject Classification}


\title[A TS algorithm for control problems with state constraints]{A tree structure algorithm for optimal control problems with state constraints}

\author[A. Alla]{Alessandro Alla}
\author[M. Falcone]{Maurizio Falcone}
\author[L. Saluzzi]{Luca Saluzzi}

\begin{abstract}
We present a tree structure algorithm for optimal control problems with state constraints. We prove a convergence result for a discrete time approximation of the value function based on a novel formulation of the constrained problem. Then the  Dynamic Programming approach is developed by a discretization in time leading to a tree structure in space derived by the controlled dynamics, in this construction the state constraints  are taken into account to cut several branches of the tree. Moreover, an additional pruning allows for the reduction of  the tree complexity as for the case without state constraints. Since the method does not use an a priori space grid, no interpolation is needed for the reconstruction of the value function and the accuracy essentially relies on the time step $h$. These features permit a reduction in CPU time and in memory allocations. The synthesis of optimal feedback controls is based  on the values on the tree and an interpolation on the values obtained  on the tree will be necessary if a different discretization in the control space is adopted, e.g.  to improve the accuracy of the method in the reconstruction of the optimal trajectories. Several examples show how this algorithm can be applied to problems in low dimension and compare it to a classical DP method on a grid.

\end{abstract}

\maketitle

\section{Introduction}

We deal with the following optimal control problem with state constraints.
\par
Let $\Omega$ be an open bounded subset of $\R^d$,
we consider the following system of controlled differential equations
\begin{equation}
\label{eq:dyn}
\begin{cases}
\dot y(s) = f(y(s), u(s)) \qquad s\geq 0   \cr
y(0) = x   \cr
\end{cases}
\end{equation}
Here $x\in \cGo$ and the {\it control} $u (t)$ belongs to the set of
admissible control functions $\U$, typically the set of measurable control functions with values in $U$, 
 a compact subset of $\R^m$. We impose a state constraint on
\eqref{eq:dyn} requiring that the state remains in $\cGo$ for all $t \geq 0$. As a consequence, we will consider admissible (with respect to the
state constraint) only control functions guaranteeing that the corresponding
trajectory never leaves  $\cGo$. We will denote by $\Ux$ this subset of
$\U$, then for\ any  $ x \in \cGo$
\begin{equation}
 \Ux= \{ u (\cdot) \in \U : y_x(t) \in \cGo ,\; \forall t \geq 0 \}
\end{equation}
where $y_x$ denotes the solution trajectory starting at $x$.\\
Given a cost functional $J(x,u)$, the problem is to determine the {\it
value function}  
\begin{equation}
v(x)=\inf_{u \in \Ux} {J(x,u)} , 
\end{equation}
and possibly an {\it optimal control} (at least approximate).  We
will use the notion of viscosity solution to the Hamilton-Jacobi-Bellman equation,
introduced by Crandall and Lions in \cite{CL83} (see also \cite{L82}), and in particular its extension to the notion of constrained
viscosity solution given by Soner \cite{S86} in order to treat
problems with state constraints. This definition combines the standard
definition on $\Omega$ with an appropriate inequality to be satisfied on
$\fGo$ (see also \cite{CDL90} for further developments of this notion).
This condition can be applied to other hamiltonians coming from various optimal control problems. \par

In the first part of this paper we will consider for simplicity only convex constraints for the infinite horizon problem. As we will see later, similar arguments can be applied to other control problems and non convex constraints although our result does not cover this case.
Dealing with the infinite horizon problem, Soner has shown that,
whenever the value function is continuous, it is the unique constrained
viscosity solution of the following Hamilton-Jacobi-Bellman equation 
\begin{equation}\label{HJBinf}
\lambda v(x) = \inf_{u\in U} \{f(x,u) \cdot \nabla v(x) + \ell(x,u)\}
\end{equation}
where $\lambda$ is a positive real parameter, the {\it discount rate}.\par
We should also mention that several results have been obtained for
the existence of trajectories of (1.1) satisfying the state constraints
(the so called {\it viable trajectories}) using the theory of
multivalued differential inclusions (see Aubin--Cellina \cite{AC84}).
Essentially, we know that a viable solution exists if for any $x \in
\cGo$ there exists at least one control such that the corresponding
velocity $f(x,u)$ belongs to the tangent cone to $\Omega$ at $x$ (see
Section 2 for a precise result in the convex case due to Haddad \cite{H81}). We recall that  several extensions have been proposed
for more general constraints using appropriate definitions of tangent cones
(see \cite{A91} for an extensive presentation of this theory).
These results gives necessary and sufficient
conditions for the existence of viable trajectories so that one can determine the minimum set of assumptions guaranteeing that
the optimal control problem can have a solution.\par
Several papers have been written on optimal control problems with state constraints starting from 
the seminal paper  \cite{S86}. We can mention the interesting contributions by Ishii-Koike \cite{IK96}, Bokanowski-Zidani and co-authors \cite{BFZ11, ABZ13} and Motta \cite{M95} for different ways to deal with state constraints still having a well posed problem. We also mention the recent contribution by Kim, Tran and Tu \cite{KTT20} dealing with constrained problems on nested domains.\\
From the point of view of the numerical approximation a classical grid approach has been developed 
by Camilli-Falcone \cite{CF96}   and Bokanowski-Forcadel-Zidani \cite{BFZ11}.\\
In this respect they represent an extension to constrained problems of the numerical approximation developed by Capuzzo Dolcetta \cite{CD83},
Falcone \cite{F87} (see also the survey paper \cite{CDF89} and the book \cite{FF13} for other numerical schemes related to optimal control problems via the Dynamic Programming approach). We end this short presentation mentioning that also viability tools have been applied to construct numerical methods for optimal control problems with state constraints, see e.g. \cite{CQSP97}.\\
Although convergence results are available for every dimension, numerical methods based on fixed space grids 
are difficult to apply for high-dimensional problems since they suffer for the well known 'curse of dimensionality'.
This is why a renewed effort has been made in recent years to find other methods which can tackle high-dimensional optimal control problems.
A list of references for other approaches dealing with high-dimensional problems is presented and discussed in \cite{AFS18}.

In the first part of this paper we propose a novel formulation of the time discrete infinite horizon problem that is close to the formulation presented in \cite{IK96} for the continuous problem and we prove a convergence result for the value function 
for a convex constraint. The proof is based on a mixture of tools coming from  multivalued analysis and viscosity solutions. As we said, we want to develop a fast approximation scheme for the value function using the characterization in terms of the Hamilton-Jacobi-Bellman equation. To this end we will also use some tools of the viability theory to establish a precise convergence result (see Sections 2 and 3). The scheme is build having in mind a "heuristic" representation of the value function which comes out coupling the viability results with standard dynamic
programming arguments.
 Although we present our convergence result for the infinite horizon problem focusing on the  treatment of boundary conditions for the stationary problem, similar arguments can be applied  also to other optimal control problems such as the finite horizon and the optimal stopping problem (see Remark \ref{genpro}).\\
The second part of the paper is devoted to the construction of an efficient algorithm for a  time discrete approximation of the value function that avoids the construction of a fixed grid in space and allows to apply the dynamic  programming principle on a Tree Structure (TS), the main results on this approach have been presented \cite{AFS18, AFS18b, AS19}. Our contribution here is the extension of the TS Algorithm (TSA)  to problems with state constraints and the feedback reconstruction using scattered data interpolation.

The outline of our paper is the following.\par\noindent
In Section 2 we introduce our basic assumptions and state some
previous results about the characterization of
the value function in terms of the Hamilton-Jacobi-Bellman equation. We present some
results in the viability theory that are useful for the problem at
hand and  discuss a different way to write the equation. We continue introducing  our time discretization and prove some properties
of the discrete value function $v_h$ showing  that the discretized equation \eqref{HJBh}
has a unique solution $v_h$.  We establish our main convergence result for the infinite horizon problem in Section 3 proving that $v_h$ converges to the value $v$ uniformly on the constraint $\cGo$, provided the state constraint $\Omega$ is convex. In Section 4 we introduce the TSA for the finite horizon problem with state constraints and discuss some of its features. Finally,  the last  section is devoted to numerical experiments where we show the TSA is faster than the classical grid approximation.  Moreover, some of the tests show that the method can also solve problems with non convex space constraints, overcoming the limits of our convergence result.

\section{The infinite horizon problem with state constraints.} 

We will denote by $y(x,t,\overline{u} (t))$ the position at time $t$ of
the solution trajectory of (1.1) corresponding to the control
$\overline{u} \in \mathcal{U}$. Whenever this will be
possible without ambiguity we adopt the simplified notations
$y_x(t)$ or $y(t)$ instead of  $y(x,t,\overline{u} (t))$. The cost functional related to the infinite horizon
problem is given by
\begin{equation}\label{def:orinf}
J(x,u) \equiv \int_{0}^{+\infty} \ell(y(t),u(t)) e^{-\lambda t}dt ,
\end{equation}
where $\ell$ is the {\it running cost}. As we said in the introduction  
we want to minimize $J$ with respect to the controls in $\Ux$ so we need
at least the assumption that 
\begin{equation}\label{as:Unonempty}
\Ux \not = \emptyset \qquad \ {\rm for \ any \ } x \in \cGo.
\end{equation}
It is important to note that in general $v(x)$ is not continuous on
$\cGo$ even when \eqref{as:Unonempty} is satisfied. This is due to the structure of the
multivalued map  $x \rightarrow \Ux $.

Soner has shown in \cite{S86} that the value function is continuous (and then 
uniformly continuous) on $\cGo$ if the following boundary condition on
the vectorfield is satisfied
\begin{equation}\label{def:sonerc}
\exists\gamma>0\;: \forall x\in\fGo \; \exists u\in U \; {\rm such\ that \ } f(x,u)\cdot\eta(x)\leq-\gamma<0
\end{equation}
where $\eta(x)$ is the outward normal to $\Omega$ at the point $x$.\par
 
We will make the following {\em assumptions}:\\
{\em A0.} $\Omega$  is a bounded, open  convex  subset  of  $\R^d$; \\
{\em A1. }$U\subset\R^d$,  compact;\\
{\em A2.} $f: \R^d\times U\longrightarrow \R^d$, is continuous   and  $\sup\limits_{u\in U}|f(x,u)-f(y,u)|\leq L_f|x-y|,$\\
{\em A3.} $ \ell : \R^d\times U\longrightarrow \R$  is continuous and  $\sup\limits_{u\in U}|\ell(x,u)-\ell(y,u)|\leq L_\ell|x-y|$.\\

Clearly, there exist two positive constants $M_\ell$, $M_f$ such that
\begin{equation}\label{bounds}
 \sup\limits_{u\in U}|f(x,u)|\leq M_f    \hbox{ and } \sup\limits_{u\in U}|\ell(x,u)|\leq M_\ell                      
 \end{equation}
for any $x \in \cGo$. Notice that under the above assumptions the value
function is bounded in $\cGo$ by $M_\ell/\lambda$ as can be easily checked.

Using the Dynamic Programming Principle, Soner has shown that $v$ is the unique viscosity solution of  \eqref{HJBinf}. This means that $v$
satisfies 
\begin{equation}\label{eq:cvis1}
H(x,u(x),\nabla u(x))\le 0\qquad {\rm for} \; x\in \Omega
\end{equation}
\begin{equation}\label{eq:cvis2}
H(x,u(x),\nabla u(x))\ge 0\qquad {\rm for} \; x\in \cGo
\end{equation}
where 
\begin{equation}\label{eq:ham}
H(x,u(x),\nabla u(x))\equiv \lambda u(x)+\max_{u\in U} \{-f(x,a)\cdot \nabla u(x) -\ell(x,a)\} 
\end{equation}
and the above inequalities should be understood in the viscosity sense (see \cite{S86} for the precise definition). A function satisfying \eqref{eq:cvis1} (respectively
\eqref{eq:cvis2}) is be called a {\it constrained viscosity subsolution} (respectively {\it supersolution}) of $H(x,u(x),\nabla u(x))=0$. \par

\begin{teo}\label{teo0}
Let \eqref{as:Unonempty},  (A0) -(A3) be satisfied and let us assume that $v
\in C(\cGo)$. Then, $v$ is the unique viscosity solution of  \eqref{HJBinf} on $\cGo$.
\end{teo}

\begin{rmk}{Necessary and sufficient conditions.}\\
Condition \eqref{def:sonerc} is known to be only a sufficient condition for the
existence of trajectories living in $\cGo$. However, necessary and
sufficient condition for the existence of solutions in $\cGo$ have been
extensively studied in {\it viability theory} (see \cite{A91}).\\
Let  $\Omega$ be an  open convex subset of $\R^d$. A trajectory is called  {\em viable} when 
\begin{equation}\label{viab}
y(t) \in \cGo, \qquad \forall t\geq 0.
\end{equation}
\end{rmk}

Let $F : \cGo \rightarrow \R^d$ be a multivalued map which is lower semicontinuous and has compact
convex images (we refer to \cite{AC84} for the theory and the definitions related to multivalued maps). Let us define the
tangent cone to a compact convex set $K$ at the point $x$, as 
\begin{equation} \label{def:tcone}
T_{K}(x)\equiv {\rm cl}\left(\bigcup_{h>0}{1\over h}(K-x)\right).
\end{equation}
A result due to Haddad \cite{H81} shows that the condition 
\begin{equation}\label{(FiT)}
F(x) \cap T_{\cGo} \neq\emptyset , \qquad \forall x\in \cGo ,
\end{equation}
is necessary and sufficient to have viable trajectories in $\cGo$ for the multivalued Cauchy problem
\begin{equation}
\begin{cases}
\dot y(t) \in F(y(t)) \qquad t\geq 0 \,,   \cr
y(0) = x\in \cGo \, .
\end{cases}
\end{equation}
This result has been also extended to more general sets (also non convex)  introducing more general tangent cones (see \cite{A91} for a general presentation of the viability theory). 

\subsection{The time-discrete scheme for the constrained problem} 
In order to build a discretization of  \eqref{HJBinf} we start using the standard
discretization in time of \eqref{eq:dyn}, \eqref{def:orinf}.  
We fix a positive parameter $h$, the time step, and consider the
following approximation scheme for (1.1)  and (2.1)
\begin{equation}
\begin{cases}
y_{n+1}=y_n + hf(y_n,u_n), \qquad n\in\mathbb{N}   \cr
y_0=x     
\end{cases}
\end{equation}
\begin{equation} \label{def:jh}
J^h(x,\{u_n\} )=h\sum_{n=0}^{+\infty} f(y_n,u_n )\beta ^k,
\end{equation}
where $x \in \cGo$, $u_n\in U$ and $\beta \equiv1-\lambda h$.\\
For every $x \in \cGo$ the  corresponding value function is
\begin{equation} \label{def:vh}
v_h(x)=\inf_{\{u_n\} \in \U^h_x } J^h(x, \{u_n\} ),
\end{equation}
where 
\begin{equation}
\U^h_x = \{ \{u_n\} : u_n \in U \hbox{ and } y_n \in \Omega,  \;\; \forall n\in\mathbb{N}\}
\end{equation}

The above definition is  meaningful only provided there exists a step $h$ such that $\U^h_x \neq \emptyset $. We look for
conditions guaranteeing the existence of viable discrete trajectories. Let us introduce the multivalued map
\begin{equation} \label{def:Uhx}
 U_h(x) \equiv \{ u \in U :  x+hf(x,u) \in \Omega \}.
 \end{equation}
representing  the subset of admissible (i.e. satisfying the constraint) controls for the discrete dynamics.
Clearly $\{u_n\} \in \U^h_x$ if and only if $u_n \in U_h(y_n)$ for any $n\in \mathbb{N}$. Due to the regularity assumptions on $f$, $U_h(x)$ is
open and is bounded since  is always contained in $U$.

\vspace{0.3cm}
\begin{rmk}\label{rem32}
Note that 
\begin{equation} \label{}
\hbox{if } u \in U_h(x), \hbox{ then } f(x,u) \in {\rm int} \left( T_{\cGo} (x)\right)
\end{equation}
where ${\rm int}\left( T_{\cGo} (x)\right)$ is the interior of the tangent cone to $\cGo$ at $x$, i.e.
\begin{equation}
 {\rm int}\left( T_{\cGo} (x)\right) = \bigcup_{h>0}\frac{1}{h} (\Omega-x).
\label{3.7}
\end{equation}
 In fact, if $u \in U_h(x)$, then $x+f(x,u) \in \Omega$, which implies  $f(x,u) \in\frac{1}{h}(\Omega-x) \subset{\rm int}\left( T_{\cGo} (x)\right)$. Note that
${\rm int}\left( T_{\cGo} (x)\right)$ is not empty since $\Omega \neq \emptyset$.

The dependence of $U_h(x)$ from $h$ is such that 
\begin{equation}\label{3.8}
U_h(x) \subset U_t(x) \qquad \forall t \in (0,h], \forall x \in \cGo .
\end{equation}
In fact, if $u \in U_h(x)$ then $x+f(x,u) \in \Omega$ and \eqref{3.8} follows by the convexity of $\cGo$.
\end{rmk}

The following proposition gives necessary  and sufficient conditions  for the existence of a time step  
$h$, such that $U_h(x) \neq \emptyset$ for any $x\in \cGo$ and therefore guarantees $\U^h_x \neq \emptyset$.

\begin{prp} \label{prop1}
Let $\Omega$ be an open bounded convex subset of $\R^d$.
Assume that $f : \cGo \times U \rightarrow \R^d$ is continuous.
Then,  there exists $h>0$ such that 
\begin{equation}\label{(3.9)}
 U_h(x) \neq \emptyset \qquad  \hbox{ for any } x\in \cGo 
 \end{equation}
if and only if the following assumption holds,
\begin{equation} \label{as:A4}
\forall x \in \fGo, \; \exists u \in U : f(x,u) \in {\rm int}\left(T_{\cGo} (x)\right).
\end{equation}
\end{prp}

\begin{proof}
If such an $h>0$ exists, \eqref{as:A4} is satisfied by Remark \ref{rem32}.\\
Now let us consider an $x \in \fGo$ and let $u=u(x) \in U$ be a control
satisfying \eqref{as:A4}. Since  $f(x,u) \in {\rm int}\left(T_{\cGo} (x)\right)$ and $f$ is bounded there exists an
$h_{x,u}>0$ such that 
\begin{equation} \label{rel1}
x+h_{x,u}f(x,u) \in \Omega .
\end{equation}
Moreover, \eqref{rel1} is also valid for every positive $h\le h_{x,u}$ by the convexity of $\Omega$.
Since $f$ is bounded, \eqref{rel1} is satisfied for any $x \in \Omega$ and $h_{x,u}$ will not depend on $x$. By the continuity of $f$  there will be  an $h$ and a neighbourhood $I(x)$ of $x$ such that
$$
\forall y \in I(x) \cap \cGo, \qquad y+h f(y,u) \in \Omega
$$
at least for $u=u(x) \in U$.\par\noindent 
We define
\[
\mathcal{O}_h \equiv \{x\in \cGo \; \vert \; \exists u \in U : x+hf(x,u) \in \Omega  \hbox{ for an } h>0\}. 
\]
Note that when $x\in \Omega$ all the directions are allowed provided $h$ is sufficiently small and the restrictions apply only for $x\in\partial \Omega$.
The family $\mathcal{O}_h$ is an open covering of $\cGo$ from which we can
extract a finite covering $\{ \mathcal{O}_{h_j} \}_{j=1,\ldots,p}$. We will
have then  $U_h(x) \neq \emptyset$ for any  $x\in \cGo$ setting $h={\displaystyle \min_j} \{h_j\}$.  
\end{proof}

\begin{cor}\label{cor1}
Under the same assumptions of Proposition \ref{prop1} there exists  $h>0$ such that
\begin{equation}
U_t(x) \neq \emptyset \qquad \forall t\in (0,h], \forall x \in \cGo .
\label{(3.11)}
\end{equation}
\end{cor}
Let us remark  that condition \eqref{as:A4} is more general than the boundary condition 
\eqref{def:sonerc} since does not require the regularity of  $\partial\Omega$.  In fact for  a closed convex subset $K$, we can define the normal cone $N_{K}(x) $ to $K$ at as
\begin{equation}\label{def:ncone}
N_{K}(x) \equiv\{y\in \R^d : \langle y,z\rangle \le 0 ,\quad \forall z\in T_{K}(x) \}
\end{equation}
When $x\in\Omega$ the tangent cone will be the whole space $\R^d$ and the normal cone will be empty. For $x\in\partial K$ these are real convex cones.
Now assume that $\cGo$ has a regular boundary, the tangent cone is an hyperplane
and the normal cone is reduced to  $\lambda \eta (x)$, $\lambda >0$.  Then \eqref{def:sonerc}
implies that 
\begin{align}
 &{\forall x\in\fGo \ \exists\  u=u(x) \in U \hbox{ such  that}}  \nonumber\\
&\langle f(x,u),v\rangle =\langle f(x,u),\lambda \eta (x)\rangle \leq -\lambda \gamma <0 \qquad \forall v \in N_{\cGo}(x) 
\end{align}
hence
$$
f(x,u) \in {\rm int}\left(T_{\cGo}(x)\right).
$$

In the sequel we will use condition \eqref{as:A4} instead of \eqref{def:sonerc}. \par
The proof of the following result can be obtained by standard arguments so it will not be given here (see \cite{BCD97} for details).

\begin{prp} \label{prop2}
\begin{equation}\label{DDPP}
v_h(x)=\inf_{\{u_n\} \in \U^h_x} \left( h\sum^{p-1}_{k=0} \ell(y_k,u_k)\beta^k+\beta^p v_h(y_p)\right),
\end{equation}
for any $x\in \cGo$ and $p\geq 1$, where $y_k$ is the trajectory with the the con sequence $\{u_j\}_{j=0}^k$.
\end{prp}

We will refer to \eqref{DDPP} as the Discrete Dynamic Programming Principle
(DDPP). For $p=1$,it gives the following discrete version of $\eqref{HJBinf}$ 

\begin{equation}\label{HJBh}
v(x)=\inf_{u \in U_h(x)} \{\beta v(x+hf(x,u))+h\ell(x,u) \}, \qquad x\in \cGo.
\end{equation}
Note that for the constrained problem the infimum is taken on the variable control set $U_h(x)$. In the next section  we will
see how to handle this dependency.
\begin{teo}\label {teo2}
Let $\lambda > L_f$. Then, for any $h \in (0,\frac{1} {\lambda}]$ there exists a unique solution $v_h \in C(\cGo)$ of $\eqref{HJBh}$.
Moreover, the following estimates hold true: 
\begin{equation} \label{est1}
 \omega_{v_h}(\delta) \leq \frac {L_\ell} {\lambda -L_f} \delta , \qquad \delta >0
\end{equation}
\begin{equation}\label{est2}
\Vert v_h \Vert _{\infty} \leq {M_\ell \over \lambda}.
\end{equation}
where $\omega_{v_h}$ is the modulus of continuity of $v_h$.
\end{teo}

\begin{proof}
The solution of \eqref{HJBh} is the fixed point of the operator $T$
\label{ (3.16)}
\begin{equation}\label{def:T}
Tv(x)=\inf_{u \in U_h(x)} \{ \beta v(x+hf(x,u) )+h\ell(x,u) \}, \qquad x \in \cGo.
\end{equation}
Let $u,v \in L^{\infty}(\cGo)$ and $x \in \cGo$. By (3.16)
for any $\varepsilon >0$, there exists $u^{\varepsilon}=u^{\varepsilon}(x,v) \in U_h(x)$ such that
\begin{equation}\label{ineq1}
Tv(x)+\varepsilon \geq  \beta v(x+f(x,u^\varepsilon) )+h\ell(x,u^\varepsilon) ,
\end{equation}
then
\begin{eqnarray}
Tu(x)-Tv(x)&\leq \beta [u(x+f(x,u^\varepsilon) )-v(x+f(x,u^\varepsilon) )]+ \\
&+h[\ell(x,u^\varepsilon)-\ell(x,u^\varepsilon)]+\varepsilon \leq \beta \Vert u-v\Vert _{\infty} +\varepsilon,  \nonumber
\end{eqnarray}
which implies 
\[
Tu(x)-Tv(x) \leq \beta \Vert u-v\Vert _{\infty} ;
\]
Reversing the role of $u$ e $v$ we get
\begin{equation}\label{contrac}
\Vert Tu-Tv\Vert _{\infty} \leq \beta \Vert u-v\Vert _{\infty}.
\end{equation}
Note that if  $v\in L^\infty(\overline \Omega)$ is such that $\Vert v \Vert _{\infty} \leq M$, we have 
\[
|Tv(x)| \leq \beta \Vert v \Vert _{\infty}+hM_\ell \leq \beta M+hM_\ell;
\]
Then, recalling the definition of $\beta$,  $\Vert v \Vert _{\infty} \leq {M_\ell \over \lambda}$ implies
\begin{equation}
\Vert Tv\Vert _{\infty} \leq {M_\ell \over \lambda}.
\end{equation}
We can conclude that, for any $h \in (0,\frac{1}{\lambda}]$, $T$ is a contraction mapping in $L^{\infty}(\cGo)$ so that there will be a
unique bounded solution $v_h$ of $\eqref{HJBinf}$.\par

Now we prove that $v_h\in C(\cGo)$. We show first that if $v \in C(\cGo)$ then $Tv \in C(\cGo)$. Let $x \in \cGo$, for any
$\varepsilon >0$ there exists $u^\varepsilon =u^\varepsilon (x,v) \in U_h(x)$ which satisfies \eqref{ineq1}. Since $\Omega$ is open and  $f$ is continuous, there will be a neighbourhood  $I(x)$ of $x$ such that 
\[ 
\forall y \in I(x) \cap \cGo, \qquad y+hf(y,u^\varepsilon) \in \Omega, 
\]
then $u^\varepsilon \in U^h(y)$ and we have
\begin{equation}\label{ineq1b}
Tv(y) \leq \{ \beta v(y+hf(y,u^\varepsilon) )+hf(y,u^\varepsilon) \}.
\end{equation}
By \eqref{ineq1} and  \eqref{ineq1b} we get
\begin{eqnarray}\label{}
&Tv(y)-Tv(x) \leq \beta [v(y+hf(y,u^\varepsilon) )-v(x+hf(x,u^\varepsilon) )]+h[\ell(y,u^\varepsilon)-\ell(x,u^\varepsilon)]+\varepsilon  \nonumber \\
&\leq \beta \omega_v((1+hL_f)|x-y|)+hL_\ell|x-y|+\varepsilon   \nonumber
\end{eqnarray}
where 
\[
|y+hf(x,u^\varepsilon)-x-hf(x,u^\varepsilon)| \leq (1+hL_f)|x-y|
\]
By the arbitrariness of $\varepsilon$, we conclude 
\[
Tv(y)-Tv(x) \leq \beta \omega_v((1+hL_f)|x-y|)+hL_\ell|x-y|.
\]
Since $x$ and $y$ are arbitrary, we can determine $\delta >0$ such
that 
\begin{equation}\label{ineq3}
|Tv(y)-Tv(x)| \leq \beta \omega_v((1+hL_f)|x-y|)+hL_\ell|x-y|
\end{equation}
whenever $|x-y| \leq \delta$. By \eqref{ineq3} we get
$$
\omega_{Tv}(\delta) \leq \beta \omega_v((1+hL_f)\delta )+hL_\ell \delta
$$
and by the uniform continuity of $v$
$$
\lim_{\delta \to 0^+} \omega_{Tv}(\delta) = 0,
$$
then  $Tv \in C(\cGo)$.\par
Since $\lambda>L_f$, the constant $C_h = {hL_\ell \over 1-\beta (1+hL_f)}$ is strictly positive and one can easily check that 
$$
\omega_{Tv_0}(\delta) \leq C_h\delta ,
$$
for any $v_0 \in W^{1,\infty}(\cGo)$ such that $\omega_{v_0}(\delta) \leq
C_h\delta$. Then the recursion sequence 
$$
v_1=Tv_0, \qquad v_n=Tv_{n-1} \qquad n=2,3,\ldots
$$
starting at a $v_0$ such that $\Vert v_0\Vert_{\infty} \leq {M_\ell \over
\lambda}$ and $\omega_{v_0}(\delta)  \leq C_h\delta$ converges to the unique solution
$v_h \in L^\infty(\cGo)$ of \eqref{HJBh}. By (3.19) $v_h$ satisfies (3.15).  Since $C_h$ is decreasing in $h$,
we get
$$
\omega_{v_h}(\delta) \leq C_h\delta \leq \max_{h>0} {hL_\ell \over 1-\beta (1+hL_f)}\delta ={L_\ell \over \lambda -L_f} \delta,
$$
and we can conclude  the proof  of the theorem.
\end{proof}

\section{A convergence result}

The main result of this section is that the solution $v_h$ of the discrete--time equation converges to $v$. In order to prove this convergence we need some preliminary lemmas on the regularity of $U_h(x)$ with respect to $h$.

\begin{prp}\label{prop3}
For any fixed $h>0$, the multivalued map $x \to U_h(x)$, $x \in \cGo$, is lower semicontinuous $(l.s.c.)$ in the sense of multivalued maps
\end{prp}
\begin{proof}
Let $u_x \in U_h(x)$ and $\varepsilon >0$. Recalling  the definition of l.s.c. maps ( (see \cite{AC84}), we have to show that there exists a neighborhood $I(x)$ of  $x \in \cGo$ such that 
\begin{equation}\label{rel2}
\forall y \in I(x) \qquad \exists u_y \in U^h(y) \cap (u_x + \varepsilon B).
\end{equation}
where $B$ is the unit ball in $\R^d$.
Since $\Omega$ is open and $f$ is continuous, we can determine $\delta_1 >0$ and $\varepsilon _1 >0$ such that 
\begin{equation}\label{rel3}
\forall y \in (x + \delta _1 B) \cap \cGo,\; \forall u\in (u_x + \varepsilon _1 B) \cap U, \qquad y+hf(y,u) \in \Omega,
\end{equation}
then  $u \in U^h(y)$. Then we take $\varepsilon _1 <\varepsilon$ and $\delta _1 >0$
such that \eqref{rel3}  holds and we get \eqref{rel2} setting $I(x)=x + \delta _1 B$.  
\end{proof}

\begin{teo}\label{teo3}
 Let  $x \in \Omega$ and consider the sequence of sets $\{U_{h_p}(x)\}_p$, $p \in {\bf N}$.  Let $h_p \to 0^+$ per
$p \to +\infty$, then 
\begin{equation} \label{rel4}
U \subset \underline{\rm Lim} \{U_{h_p}(x)\} \qquad {\rm for\ } p\to +\infty .
\end{equation}
\end{teo}

\begin{proof}
Let $u \in U$, we have to prove that $u \in \underline{\rm Lim}
\{U_{h_p}(x)\}$, i.e. that for any $\varepsilon >0$, there exists an
index $\overline p$ such that 
\begin{equation}
 \forall p\geq \overline p,
\qquad U_{h_p}(x) \cap (u + \varepsilon B) \neq \emptyset.
\end{equation}
Since $x \in \Omega$ and $f$ is bounded, there exists $h_{x,u} >0$
such that 
$$
x+h_{x,u}f(x,u) \in \Omega.
$$
By a compactness argument we can choose $h_{x,u}$ independently of $u$. The continuity of $f$ then implies that there exists
$\delta >0$ such that 
\begin{equation}\label{rel5}
\forall u^{\prime} \in (u + \delta B), \qquad x+h_x f(x,u^{\prime}) \in \Omega.
\end{equation}
Moreover,  there exists an index $\overline p$ such
that 
$$
\forall p\geq \overline p, \qquad 0<h_p<h_x,
$$
then by the convexity of $\Omega$ also
$$
x+h_p f(x,u^{\prime}) \in \Omega,
$$
so  $u^{\prime} \in U^{h_p}(x)$. To end the proof it suffices to
choose $\delta <\varepsilon$  such that \eqref{rel5} holds.   
\end{proof}
Using the above propositions, we can prove our main convergence result.
\begin{teo}
Let $\lambda>L_f$, then $v_h \to v$ uniformly in $\cGo$, for $h \to 0^+$.
\end{teo}
\begin{proof}
Since $v_h$ is uniformly bounded and equicontinuous, by the Ascoli--Arzel\`a theorem, there exist $h_p \to 0^+$ for $p \to +\infty$ and a function $v \in C(\cGo)$ such that 
\begin{equation}
v_{h_p} \to v \qquad {\rm per\ } p\to +\infty \qquad \hbox{ uniformly on } \cGo\,.
\end{equation}
We will show that $v$ is the constrained viscosity solution
of \eqref{HJBinf} in $\cGo$. \\
{\em a) Let us prove first that $v$ is a subsolution of \eqref{HJBinf} in $\Omega$.}\\
Let $\phi \in C^1(\cGo)$ and let $x_0 \in \omega$ be a strict local maximum point for $v-\phi$ in
$\omega$, we have then
$$ (v-\phi )(x_0) > (v-\phi )(x) \qquad \forall x \in B(x_0,r) \subset \Omega 
$$ 
for $r>0$ sufficiently small. Then,  for $p$ large enough, there exists $x_0^{h_p} \in B(x_0,r)$ such that $v_{h_p}-\phi$ has a
local maximum point at $x_0^{h_p}$ and $x_0^{h_p}$ converges to $x_0$.
Note that for any control  $u \in U^{h_p}(x_0^{h_p})$ the point
$x_0^{h_p}+h_pf(x_0^{h_p},u)$ belongs to $\Omega$, and for $p$ large enough it
belongs to $B(x_0,r)$. The above remarks imply
\begin{equation} \label{ineq4}
v_{h_p} (x_0^{h_p} )-\phi (x_0^{h_p}) \geq v_{h_p} (x_0^{h_p}+h_pf(x_0^{h_p},u) )-\phi (x_0^{h_p}+h_pf(x_0^{h_p},u)).
\end{equation} 
By \eqref{ineq4} and  \eqref{HJBh} we get 
\begin{eqnarray}
&0 =v_{h_p} (x_0^{h_p})+\sup\limits_{u \in U_{h_p}(x_0^{h_p})} \{ -(1-\lambda h_p )v_{h_p} (x_0^{h_p}+ h_p f(x_0^{h_p} ,u))-h_p \ell(x_0^{h_p} ,u) \}\nonumber\\
& \geq  \sup\limits_{u \in U_{h_p}(x_0^{h_p} )}\{ \phi (x_0^{h_p})-\phi (x_0^{h_p} + h_p f(x_0^{h_p} ,u))+\lambda h_p v_{h_p} (x_0^{h_p}+ h_p f(x_0^{h_p} ,u))-h_p \ell(x_0^{h_p} ,u) \}\nonumber
\end{eqnarray}
Since $\phi \in C^1(\cGo)$, it follows that there exists $\theta\in [0,1]$ such that by the above inequality we get
\begin{eqnarray}\label{ineq5}
0 &\geq \sup\limits_{u \in U_{h_p}(x_0^{h_p})} \biggl \{ -\sum_{i=1}^d {\frac{\partial \phi}{ \partial x_i}}(x_0^{h_p}+ \theta h_p f(x_0^{h_p}, u)) f_i(x_0^{h_p},u)+  \\
&\lambda v_{h_p} (x_0^{h_p}+h_p f(x_0^{h_p}, u)- \ell(x_0^{h_p} ,u) \biggr\}. \nonumber
\end{eqnarray}
Let $\overline p >0$ be such that $0<h_{\overline p} <h$ and $U_h(x) \neq \emptyset$. We can choose $p^\prime >0$ such that
for any $ p \geq p^\prime , \, 0<h_p<h_{\overline p}$
\begin{equation}\label{rel6}
\emptyset \neq U_{h_{\overline p}}(x) \subset U_{h_p}(x), \qquad \forall x \in \cGo, \hbox{ for }p\geq p^\prime,
\end{equation}
and by \eqref{ineq5} we have
\begin{eqnarray}
0 &\geq \sup_{u \in U_{h_p}(x_0^{h_p})} \biggl\{ -\sum_{i=1}^d {\partial \phi \over 
\partial x_i}(x_0^{h_p}+\theta h_pf(x_0^{h_p},u)) f_i(x_0^{h_p},u)+         \\
&+\lambda v_{h_p} (x_0^{h_p}+h_p f(x_0^{h_p},u)- \ell(x_0^{h_p} ,u) \biggr\} . \qquad \quad {\rm for \ }p\geq p^\prime   
\end{eqnarray}
Let $x\in \cGo$ and $u \in U_{h_{\overline p}}(x)$, we define the real function $W(x,u)$, 
\begin{equation}\label{}
W(x,u)\equiv \biggl\{ -\sum_{i=1}^d {\partial \phi \over \partial x_i}
(x+\theta h_pf(x,u)) f_i(x,u) + \lambda v_{h_p} (x+h_p f(x,u))- \ell(x,u) \biggr\}  
\end{equation}
where $\theta \in [0,1]$ (note that $W$ is continuous in both
variables). Let us define 
\begin{equation}
W(x)\equiv\sup_{u \in U_{h_{\overline p}}(x)} W(x,u).
\end{equation}
By Proposition \ref{prop3} $U_{h_{\overline p}}(\cdot )$ is l.s.c. at $x_0$,
then by a standard result on multivalued map (see \cite{AC84}) $W$ is
l.s.c. at $x_0$. Since 
\begin{equation}\label{rel7}
W(x_0,u) \to \{ -\nabla \phi (x_0) \cdot f(x_0,u)+\lambda v(x_0)-\ell(x_0,u) \} \quad \hbox{  for }p\to+\infty,
\end{equation}
and $x_0^{h_p}$ converges to $x_0$, for any $\varepsilon >0$, there exists $p^{\prime\prime}>0$ such that $\forall p
\geq \max \{p^\prime, p^{\prime\prime} \}$ (4.11) and
(4.12) hold true. By the lower semicontinuity of $W$ and the
arbitrariness of $\varepsilon$ we get
\begin{equation}
0 \geq \sup_{u \in U_{h_{\overline p}}(x_0)} \{  -\nabla \phi (x_0) \cdot  f(x_0,u)+\lambda v(x_0)-\ell(x_0,u) \}.
\label{ineq13}
\end{equation}
The inequality \eqref{ineq13} is verified for any $\overline p >0$.
We show that
$$
0 \geq \sup_{u \in \underline{\rm Lim} \{ U_{h_{\overline p}}(x_0) \} } G(u) 
\qquad {\rm for\ }\overline p \to
+\infty,
$$
where
$$
G(u)\equiv -\nabla \phi (x_0) \cdot f(x_0,u)+\lambda v(x_0)-\ell(x_0,u).
$$
It suffices to prove that
$$
0\geq G(u) \qquad \hbox{  for any } u \in \underline{\rm Lim} \{
U_{h_{\overline p}}(x_0) \}. $$
In fact, for any $u \in \underline{\rm Lim} \{ U_{h_{\overline p}}(x_0) \}$, 
we can find a sequence $\{u^{h_{\overline p}} \}_{\overline p}$,
$u^{h_{\overline p}} \in U_{h_{\overline p}}(x_0)$ such that $u^{h_{\overline p}} \to u$ for
$\overline p \to +\infty$ and 
$$
0\geq G(u^{h_{\overline p}}),
$$
then passing to the limit for $\overline p \to +\infty$, by the continuity of $G$ we have 
$$
0\geq G(u).
$$
Proposition \ref{prop3} implies that
$$
U \subset \underline{\rm Lim} \{ U^{h_{\overline p}}(x_0) \} \qquad {\rm for\
}\overline p \to+\infty $$
so that
$$
0 \geq \sup_{u \in \underline{\rm Lim} \{ U^{h_{\overline p}}(x_0) \} } G(u) \geq \sup_{u \in  U} G(u).
$$

{\em b) Now we prove that $v$ is a viscosity supersolution of  \eqref{HJBinf} in $\cGo$.}\\ Let $\phi\in C^1(\cGo)$ and $x_0 \in \cGo$, be
a strict maximum point for $v-\phi$ in $\cGo$. We can use the same arguments that we used for  (4.9) in the first part of this
theorem (just replace $f(x_0,r)$ by $f(x_0,r) \cap \cGo$),  so we get
\label{ (4.14)}
 \begin{align}\label{ineq6}
 0 &\leq \sup_{u \in U_{h_p}(x_0^{h_p})} \biggl\{ -\sum_{i=1}^d {\partial \phi \over  \partial x_i}(x_0^{h_p}+\theta h_p f(x_0^{h_p},u)) f_i(x_0^{h_p},u)+  \\
&+\lambda v_{h_p} (x_0^{h_p}+h_p f(x_0^{h_p},u))- \ell(x_0^{h_p} ,u) \biggr\} .              \nonumber
\end{align}
where $\theta \in [0,1]$.

By \eqref{ineq6} for any $\varepsilon >0$ there exists $u_{h_p}^{\varepsilon}
\in U_{h_p}(x_0^{h_p})$ such that
\begin{align}\label{ineq7}
0 &\leq \sup_{u \in U_{h_p}(x_0^{h_p})} \biggl\{ -\sum_{i=1}^d {\partial \phi \over 
\partial x_i}(x_0^{h_p}+\theta h_p f(x_0^{h_p},u)) f_i(x_0^{h_p},u)+     \\
&\quad\quad+\lambda v_{h_p} ((x_0^{h_p}+h_p f(x_0^{h_p},u))- \ell(x_0^{h_p} ,u) \biggr\} \leq    \\
&\leq \biggl\{ -\sum_{i=1}^d {\partial\phi \over \partial x_i}f(x_0^{h_p}+\theta h_p f(x_0^{h_p},u_{h_p}^{\varepsilon})) f_i(x_0^{h_p},u_{h_p}^{\varepsilon})+               \\
&\quad\quad+\lambda v_{h_p} (x_0^{h_p}+\theta h_p f(x_0^{h_p},u_{h_p}^{\varepsilon}))- \ell(x_0^{h_p} ,u_{h_p}^{\varepsilon}) \biggr\} +\varepsilon .   
\end{align}
Since  $U_{h_p}(x_0^{h_p})$ is bounded the sequence $\{u_{h_p}^{\varepsilon}\}_p$ is also bounded. We extract a converging subsequence which we
still denote by $u_{h_p}^{\varepsilon}$. Let  $u^{\varepsilon} \in U$ be its limit for $p \to +\infty$, then  passing to the limit for $p \to +\infty$ in \eqref{ineq7} we get 
\begin{equation}
0 \leq  -\nabla \phi (x_0) \cdot f(x_0,u^{\varepsilon})+\lambda v(x_0)-\ell(x_0,u^{\varepsilon}) + \varepsilon.  
\end{equation}
and since $u^{\varepsilon} \in U$ we have 
\begin{equation}
0 \leq \sup\limits_{u \in U} \{  -\nabla \phi (x_0) \cdot f(x_0,u)+\lambda v(x_0)-\ell(x_0,u) \} + \varepsilon,
\end{equation}
Then, by the arbitrariness of $\varepsilon$ we conclude that $v$ is a viscosity supersolution of  \eqref{HJBinf} in
$\cGo$. 

Since the constrained viscosity solution is unique we conclude that $v_h$ converges  to $v$ for $h$ tending to 0.
\end{proof}

\begin{rmk}\label{genpro}
{\bf Applications to other optimal control problems.}
The above result has been proved for the infinite horizon problem, but they can be applied also to other classical control problems as we will do in the following section. The main point is in fact the boundary condition at $\partial \Omega$ that we have treated for the  infinite horizon problem.  For example, similar techniques  can be applied to solve  the {\it finite horizon} and
to the {\it optimal stopping time} problems which correspond respectively to the following cost functionals
$$ J_1(x,u) = \int_0^T \ell(y(s),u(s)) ds + g (y(T)) 
$$ 
$$ 
J_2(x,u,\tau) = \int_0^{\tau \land T} \ell(y(s),u(s))e^{-\lambda s} ds + g(y(\tau \land T)) e^{-\lambda( \tau \land T)} . 
$$
where the stopping cost $g$ is a given bounded Lipschitz
continuous function. In fact, the optimal stopping problem can
be written as an infinite horizon problem just adding a new
control, $\hat u$, to the set of admissible controls $U$ and
defining   
\[
f(x, \hat u)= 0     \hbox{ and }  \ell(x,\hat u) = {g(x) \over \lambda}. 
\]
Clearly $g$ should be Lipschitz continuous and bounded. In the following section we will present the TSA for the finite horizon problem.
\end{rmk}


\section{The finite horizon optimal control problem with state constraints}

In this section we will sketch the essential features of the dynamic programming approach for {\it finite horizon} control problems. Let the system be driven by
\begin{equation}\label{eq}
\left\{ \begin{array}{l}
\dot{y}(s)=f(y(s),u(s),s), \;\; s\in(t,T],\\
y(t)=x\in\R^d.
\end{array} \right.
\end{equation}
We will denote by $y:[t,T]\rightarrow\R^d$ the solution, by $u:[t,T]\rightarrow\R^m$ the control, by $f:\R^d\times\R^m\times[t,T]\rightarrow\R^d$ the dynamics.
We assume that there exists a unique solution for \eqref{eq} for each $u\in\mathcal{U}$.\\
 We impose the  state constraints for \eqref{eq} requiring that the state remains in a closed bounded set $\overline\Omega$ for all $t\ge0$. Our constraint is defined a rather simple way:  given an initial open set $\Omega_0$ we subtract some obstacles,  $O_j$, $j=1,\ldots,m$, that are typically defined by inequalities 
$$O_j=\{x\in \Omega_0: \varphi_j(x)\le 0\},\qquad \forall t\ge0,$$
where $\varphi(x):\R^d\rightarrow\R$ are sufficiently regular functions. We define our set $\Omega=\Omega_0\setminus \cup_{j=1}^m O_j$. Then the set of admissible controls is reduced to the following subset of $\mathcal{U}$ 
\begin{equation}
\mathcal{U}_x=\{u(\cdot)\in\mathcal{U}: y(t,u(t))\in\hbox{cl}(\Omega_0\setminus \cup_{j=1}^m O_j), \forall t\ge0\},\qquad \mbox{ for any }x\in{\overline\Omega}
\end{equation}

The cost functional for the finite horizon optimal control problem will be given by 
\begin{equation}\label{cost}
 J_{x,t}(u):=\int_t^T \ell(y(s,u),u(s),s)e^{-\lambda (s-t)}\, ds+g(y(T))e^{-\lambda (T-t)},
\end{equation}
where $\ell:\R^d\times\R^m\times [t,T]\rightarrow\R$ is the running cost. We will assume that the functions $f,\ell$ and $g$ are bounded:
 \begin{align}
 \begin{aligned}\label{Mf}
|f(x,u,s)|& \le M_f,\quad |\ell(x,u,s)| \le M_\ell,\quad |g(x)| \le M_g, \cr
&\forall\, x \in \mathbb{R}^d, u \in U \subset \mathbb{R}^m, s \in [t,T], 
\end{aligned}
\end{align}
the functions $f$ and $\ell$ are Lipschitz-continuous with respect to the first variable
\begin{align}
\begin{aligned}\label{Lf}
&|f(x,u,s)-f(y,u,s)| \le L_f |x-y|, \quad |\ell(x,u,s)-\ell(y,u,s)| \le L_\ell |x-y|,\cr
&\qquad\qquad\qquad\qquad\forall \, x,y \in \mathbb{R}^d, u \in U \subset \mathbb{R}^m, s \in [t,T], 
\end{aligned}
\end{align}
%
and the cost $g$ is also Lipschitz-continuous:
\begin{equation}
|g(x)-g(y)| \le L_g |x-y|, \quad \forall x,y \in \mathbb{R}^d.
\label{Lg}
\end{equation}
The goal is to find a state-feedback control law $u(t)=\Phi(y(t),t),$ in terms of the state equation $y(t),$ where $\Phi$ is the feedback map. To derive optimality conditions we use the well-known Dynamic Programming Principle (DPP) due to Bellman. We first define the time dependent value function for an initial condition $(x,t)\in\R^d\times [t,T]$:
\begin{equation}
v(x,t):=\inf\limits_{u\in\mathcal{U}} J_{x,t}(u)
\label{value_fun}
\end{equation}
which satisfies the DPP, i.e. for every $\tau\in [t,T]$:
\begin{equation}\label{dpp}
v(x,t)=\inf_{\color{blue} u\in\mathcal{U}_x}\left\{\int_t^\tau \ell(y(s),u(s),s) e^{-\lambda (s-t)}ds+ v(y(\tau),\tau) e^{-\lambda (\tau-t)}\right\}.
\end{equation}
Due to \eqref{dpp} we can derive the HJB for every $x\in\R^d$, $s\in [t,T)$: 
\begin{equation}\label{HJB}
\left\{
\begin{array}{ll} 
&-\dfrac{\partial v}{\partial s}(x,s) +\lambda v(x,s)+ \max\limits_{u\in U }\left\{-\ell(x, u,s)- \nabla v(x,s) \cdot f(x,u,s)\right\} = 0, \\
&v(x,T) = g(x).
\end{array}
\right.
\end{equation}
Once the value function is known, by e.g. \eqref{HJB}, then it is possible to  compute the optimal feedback control as:
\begin{equation}\label{feedback}
u^*(t):=  \argmax_{u\in U }\left\{- \nabla v(x,t) \cdot f(x,u,t)-\ell(x,u,t)\right\}. 
\end{equation}

\section{A tree structure for an optimal control problem with state constraints}

Let us start our review of the method without considering the state constraints condition. This means that $O_j=\emptyset$ and $\Omega_0 = \R^d$.
The analytical solution of Equation \eqref{HJB} is hard to find due to its nonlinearity. Here, we recall the semi-Lagrangian method on a tree structure based on the recent work \cite{AFS18}. 
Let us introduce the semi-discrete problem with a time step $h: = [(T-t)/\overline N]$ where $\overline{N}$ is the number of temporal time steps:
\begin{equation}
\left\{\begin{array}{ll}\label{SL}
V^{n}(x)=\min\limits_{u\in U}\left\{h\,\ell (x, u, t_n)+e^{-\lambda h}V^{n+1}(x+h f(x, u, t_n))\right\}, \\
\qquad\qquad\qquad\qquad\qquad\qquad\qquad\qquad\qquad\qquad\qquad\qquad n= \overline{N}-1,\dots, 0,\\
V^{\overline{N}}(x)=g(x),\quad\qquad\qquad\qquad\qquad\qquad\qquad\qquad\qquad\qquad x \in \R^d,
\end{array}\right.
\end{equation}
where $t_n=t+nh,\, t_{\overline N} = T$ and $V^n(x):=V(x, t_n).$ Note that in this section the discrete value function is denoted by $V^n(x)$ to stress the dependence on time and space. For the sake of completeness we would like to mention that a fully discrete approach is typically based on a time discretization which is projected on a fixed state-space grid of the numerical domain, see e.g. \cite{FG99}. In the current work we aim to extend the algorithm proposed in \cite{AFS18} to control problems with state constraints.

For readers convenience we now recall the tree structure algorithm. Let us assume to have a finite number of admissible controls $\{u_1,...,u_M \}$. This can be obtained discretizing the control domain $U\subset\mathbb{R}^m$ with step-size $\Delta u$.  A typical example is when  $U$ is an hypercube, discretizing in all the directions with constant step-size $\Delta u$ we get the finite set  $U^{\Delta u}=\{u_1,...,u_M \}$. To simplify the notations in the sequel we continue to denote by $U$ the discrete set of controls. Let us denote the tree by $\mathcal{T}:=\cup_{j=0}^{\overline{N}} \mathcal{T}^j,$ where each $\mathcal{T}^j$ contains the nodes of the tree correspondent to time $t_j$. The first level $\mathcal{T}^0 = \{x\}$ is clearly given by the initial condition $x$.  Starting from the initial condition $x$, we consider all the nodes obtained by the dynamics \eqref{eq} discretized using e.g. an explicit Euler scheme with different discrete controls $u_j \in U $
$$\zeta_j^1 = x+ h\, f(x,u_j,t_0),\qquad j=1,\ldots,M.$$ Therefore, we have $\mathcal{T}^1 =\{\zeta_1^1,\ldots, \zeta^1_M\}$. We note that all the nodes can be characterized by their $n-$th {\em time level}, as follows

$$\mathcal{T}^n = \{ \zeta^{n-1}_i +h f(\zeta^{n-1}_i, u_j,t_{n-1}),\, j=1, \ldots, M,\,i = 1,\ldots, M^{n-1}\}.$$

To simplify the presentation, we deal with a tree build on an Euler approximation of the dynamical system, however the algorithm can also be extended to get high-order approximations, as illustrated in \cite{AFS18b}. All the nodes of the tree can be shortly defined as

$$\mathcal{T}:= \{ \zeta_j^n,\,  j=1, \ldots, M^n,\, n=0,\ldots, \overline{N}\},$$ 

where the nodes $\zeta^n_i$ are the result of the dynamics at time $t_n$ with the controls $\{u_{j_k}\}_{k=0}^{n-1}$:
$$\zeta_{i_n}^n = \zeta_{i_{n-1}}^{n-1} + h f(\zeta_{i_{n-1}}^{n-1}, u_{j_{n-1}},t_{n-1})= x+ h \sum_{k=0}^{n-1} f(\zeta^k_{i_k}, u_{j_k},t_k), $$
with $\zeta^0 = x$, $i_k = \floor*{\dfrac{i_{k+1}}{M}}$ and $j_k\equiv i_{k+1} \mbox{mod } M$ and $\zeta_i^k \in \R^d, i=1,\ldots, M^k$. 

Although the tree structure allows to solve high dimensional problems, its construction might be expensive since $ \mathcal{T} =O( M^{\overline{N}+1}),$
where $M$ is the number of controls and $\overline{N}$ the number of time steps which might be infeasible due to the huge amount of memory allocations, if $M$ or $\overline{N}$ are too large. 
For this reason we are going to introduce the following pruning criteria:
two given nodes $\zeta^n_i$ and $\zeta^n_j$ will be merged if 
\begin{equation}\label{tol_cri}
\Vert \zeta^n_i-\zeta^n_j \Vert \le \ep, \quad \mbox{ with }i\ne j \mbox{ and } n = 0,\ldots, \overline{N}, 
\end{equation}
for a given threshold $\ep>0$. Criteria \eqref{tol_cri} will help to save a huge amount of memory. Thus, the extension to the state constraints $O_j$ case can be seen as a further pruning criteria. Indeed, together with 
\eqref{tol_cri}, we will cut off the nodes of the tree such that $\zeta_i^n\in O_j, j=1,\dots,m$.

Once the tree $\mathcal{T}$ has been built, the numerical value function $V(x,t)$ will be computed on the tree nodes in space as 
\begin{equation}\label{num:vf}
V(x,t_n)=V^n(x), \quad \forall x \in \mathcal{T}^n, 
\end{equation}
where $t_n=t+ n h$. It is now straightforward to evaluate the value function. The TSA defines a time dependent structure $\mathcal{T}^n=\{\zeta^n_j\}_{j=1}^{M^n}$ for $n=0,\ldots, \overline{N}$ and
we can obtain an approximation on the tree \eqref{dpp} as follows: 
\begin{equation}
\begin{cases}
V^{n}(\zeta^n_i)= \min\limits_{u\in U} \{e^{-\lambda h} V^{n+1}(\zeta^n_i+h f(\zeta^n_i,u,t_n)) +h \, \ell (\zeta^n_i,u,t_n) \}, \\
\qquad\qquad \qquad\qquad \qquad\qquad \qquad\qquad \qquad\qquad  \zeta^n_i \in \mathcal{T}^n\,, n = \overline{N}-1,\ldots, 0, \\
V^{\overline{N}}(\zeta^{\overline{N}}_i)= g(\zeta_i^{\overline{N}}), \qquad\qquad \qquad\qquad \qquad\qquad   \zeta_i^{\overline{N}} \in \mathcal{T}^{\overline{N}}.
\end{cases}
\label{HJBt2}
\end{equation}

We note that the minimization is computed by comparison on the discretized set of controls $U$. 

\subsection{Feedback reconstruction}

From the knowledge of the approximate value function we can obtain the synthesis feedback control. The TSA allows to store the control indices corresponding to the argmin in \eqref{HJBt2}, during the computation of the value function. Then, starting from $\zeta^0_*= x$, we follow the path of the tree to build the optimal trajectory $\{\zeta^n_*\}_{n=0}^{\overline{N}}$ in the following way
\begin{equation} \label{feed:tree}
u_{n}^{*}:=\argmin\limits_{u\in U} \left\{ e^{-\lambda h}V^{n+1}(\zeta^n_*+ h f(\zeta^n_*,u,t_n)) +h \, \ell(\zeta^n_*,u,t_n) \right\},
\end{equation}
\begin{equation*} 
\zeta^{n+1}_* \in \mathcal{T}^{n+1} \; s.t. \; \zeta^n_* \rightarrow^{u_{n}^{*}} \zeta^{n+1}_*,
\end{equation*}
for $n=0,\ldots, \overline{N}-1$, where the symbol $\rightarrow^u$ stands for the connection of two nodes by the control $u$.  We note that this is possible if we consider the same discrete control set $U$ for both HJB equation \eqref{HJBt2} and feedback reconstruction \eqref{feed:tree} as discussed in \cite{AFS18}.

In this work we are also interested to extend the feedback to a larger set of controls $\widetilde{U}$ such that $U\subset\widetilde{U}$. Therefore, the feedback control will be computed as
\begin{equation} \label{feed:tree2}
u_{n}^{*}:=\argmin\limits_{u\in \widetilde{U}} \left\{ e^{-\lambda h}I[V^{n+1}](\zeta^n_*+ h f(\zeta^n_*,u,t_n)) + h \, \ell(\zeta^n_*,u,t_n) \right\},
\end{equation}
where the $\argmin$ is computed over the new set $\widetilde{U}$ and $\zeta^n_*+ h f(\zeta^n_*,u,t_n)$ might not be a node of the tree. To this end, we need to use an interpolation operator for scattered data.
Our scattered data consists of a set of points $\{\zeta_j^n\}_{j=1}^M$ on the tree and the corresponding value function $\{V^{n+1}(\zeta_j^n)\}_{j=1}^M$, where the points have no ordering with respect to their relative locations. There are various methods for computing a polynomial interpolation on scattered data. One widely used approach uses a Delaunay triangulation of the points between the nearest point of interest and then perform a linear interpolation on the computed triangulation. Triangulation based methods are local, so they can treat  efficiently large data sets. We refer the interested reader to \cite{A02} for more details on the topic. In the numerical experiments we have used the Matlab function {\tt scatteredInterpolant}.


\section{Numerical experiments}

In this section we show our numerical results. 
The first example deals with the control of a damped harmonic oscillator in a convex constraint, where the theoretical findings obtained in the previous sections hold true. In the second example we treat a linear dynamics where we provide two different constraints: a circular channel and a labyrinth with obstacles. Finally, the last example deals with a non-linear problem, the control of the Van der Pol oscillator. The numerical simulations reported in this paper are performed on a laptop with 1CPU Intel Core i5-3, 1 GHz and 8GB RAM. The codes are written in Matlab. 

\subsection{Test 1: Damped harmonic oscillator with a convex constraint}

\begin{figure}[htbp]
\centering
\includegraphics[scale=0.4]{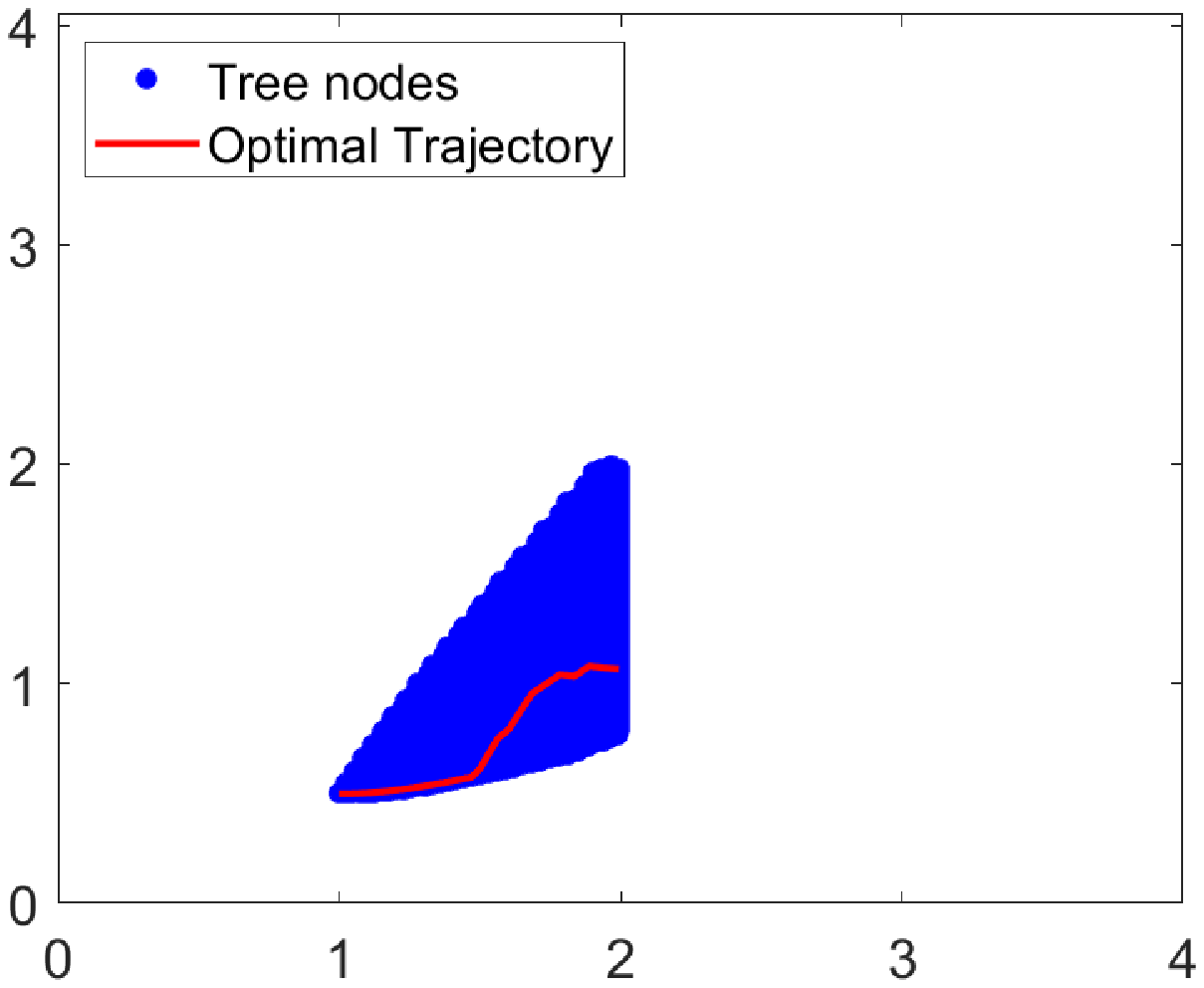}
\includegraphics[scale=0.4]{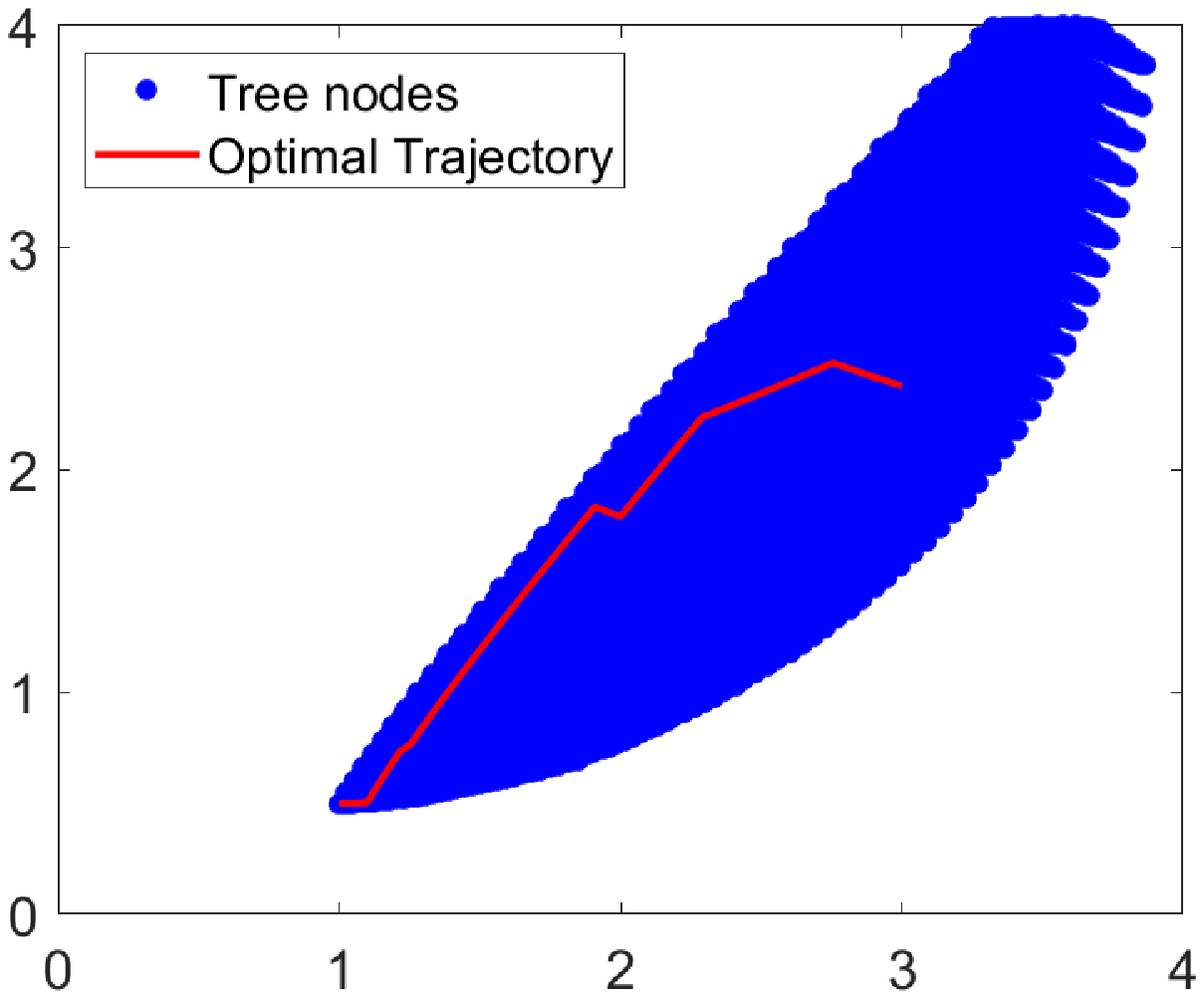}
\caption{Test 1: Optimal trajectory with constraints (left) and without constraints (right).}
\label{fig10}
\end{figure}

In the first test case we consider the damped harmonic oscillator. The dynamics in \eqref{eq} is given by
\begin{equation}\label{eq:dho}
	f(x,u)=
		  \begin{pmatrix}
         x_2 \\
         -kx_1 +u x_2
          \end{pmatrix}\quad u\in U\equiv [-1, 1].
\end{equation}
The cost functional in \eqref{cost} is: 
\begin{equation}\label{cost:2}
\ell(x,u,t) = (x_1-3)^2 \qquad g(x) = (x_1-3)^2, \qquad \lambda = 0,
\end{equation}
and aims to steer the first component of the solution to $3$. We will consider $T=1.5$ as horizon, $x=(1,0.5)$ as initial condition, $h=0.025$ and $\ep = h^2$. Furthermore, the constraint is the box $\Omega = [0,2]^2$. The optimal trajectory together with the tree nodes are shown in the left panel of Figure \ref{fig10}. We note that the solution does not reach $x_1=3$ since it is outside our constraint, but only $x_1=2$. Furthermore, the nodes of the tree do not cover the whole constraint, but only a piece. We note that the cardinality of the unconstrained tree is equal to 233739, while for the constrained tree is just 38406. For comparison on the right panel we show the solution of the unconstrained problem where we can see that the solution gets closer to the desired configuration. The number of tree of nodes is also larger since there are no restrictions and this is also reflected in the CPU of the time where we need only 7 seconds in the constraint case versus 45 second for the unconstrained. Here to build the value function and the optimal control we have used the following discrete control $\{-1,0,1\}$.

\begin{figure}[htbp]
\includegraphics[scale=0.3]{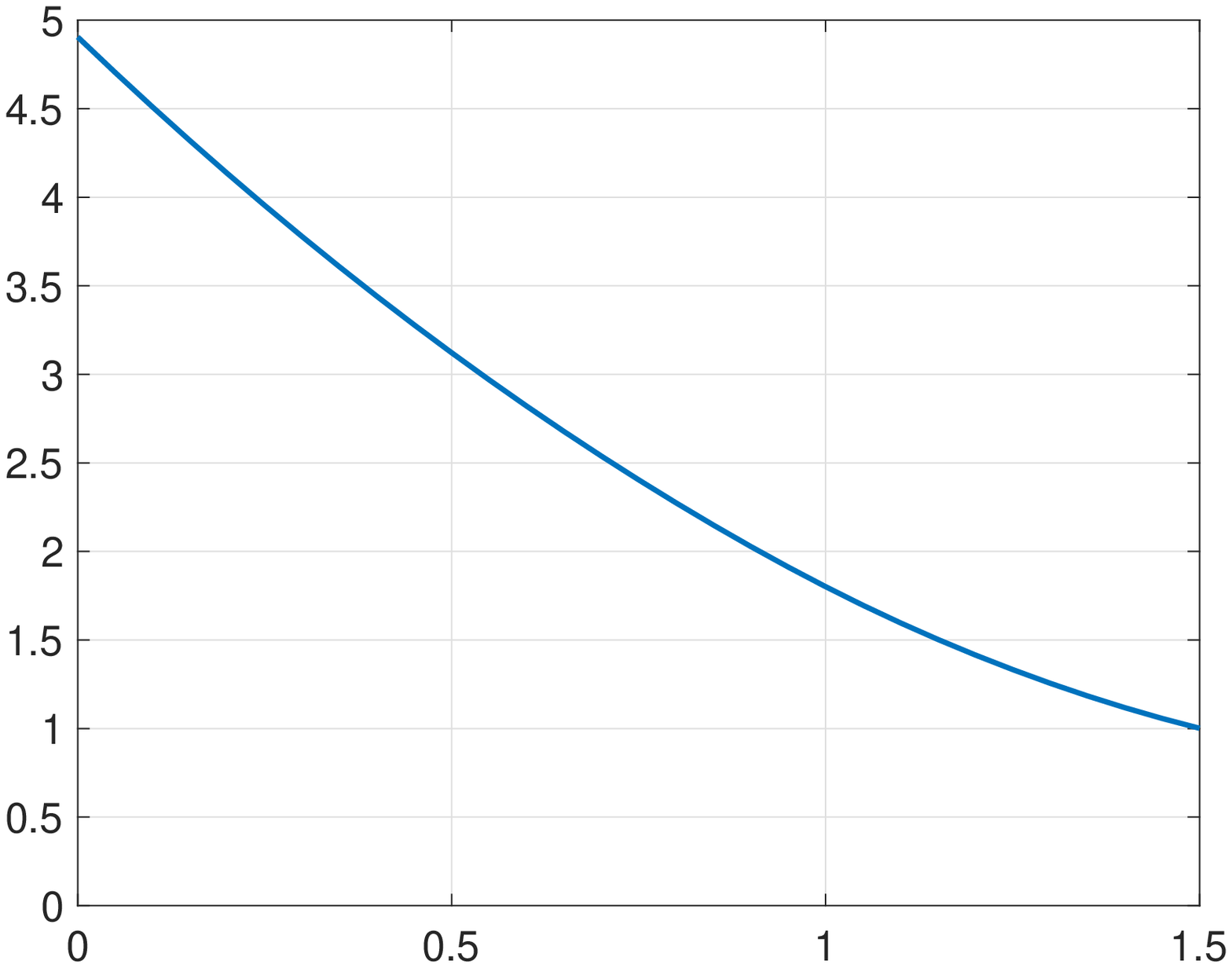}
\includegraphics[scale=0.3]{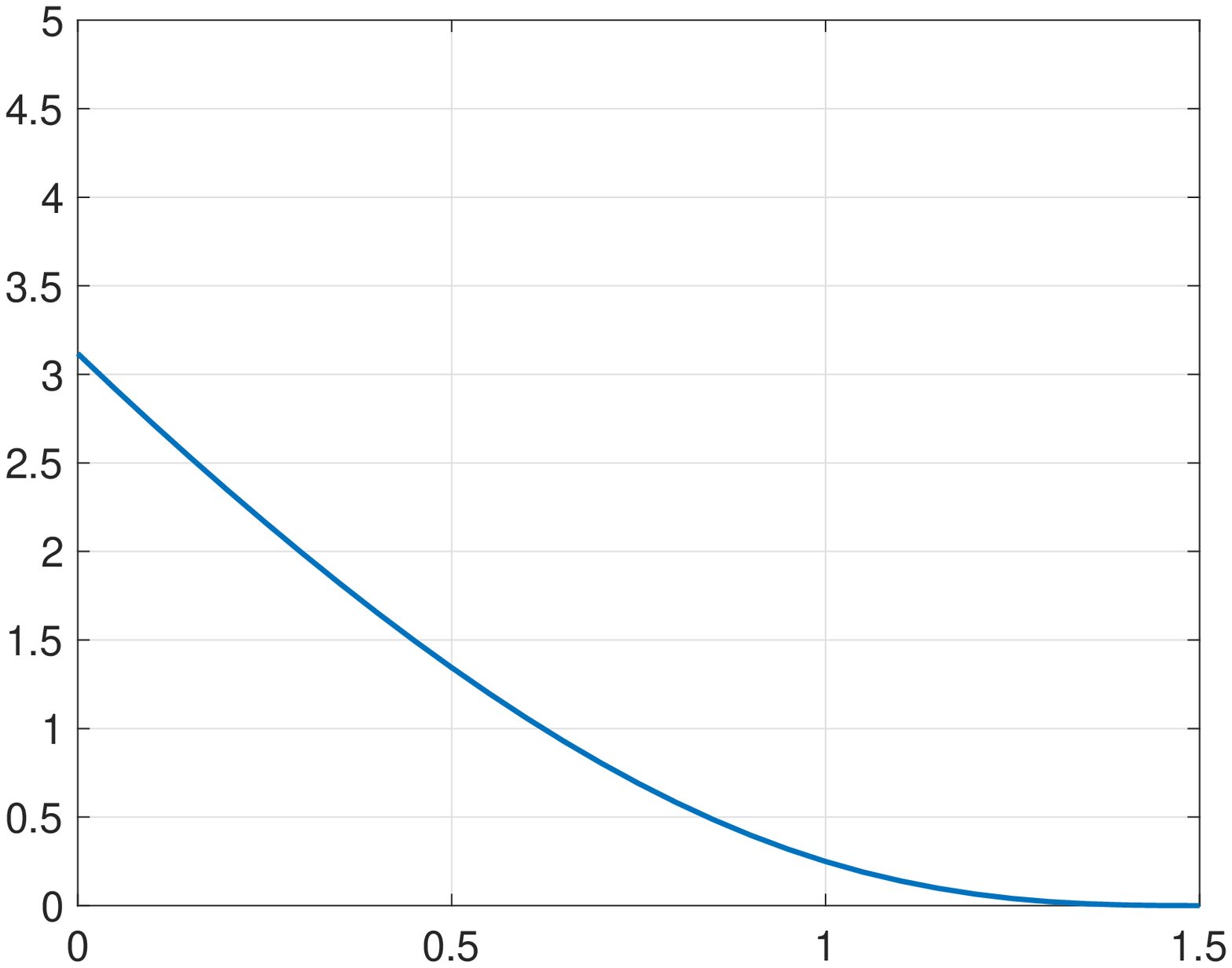}
\caption{Test 1: Cost functional with constraints (left) and without constraints (right).}
\label{fig11}
\end{figure}

The evaluation of the cost functional is shown in Figure \ref{fig11}. We can see that, as expected, the unconstrained problem (right panel) has lower value than the constraint problem (left problem). This is due to the fact that we only reach our target in the unconstrained case.

%
%
\subsection{Test 2: Simple dynamics}
In the second example, we consider the following dynamics in \eqref{eq}
\begin{equation}
	f(x,u)=
		  \begin{pmatrix}
         u_1 \\
         u_2
          \end{pmatrix},\, u \in U\equiv \partial B(0,1) \cup (0,0),
       \label{exact}   
\end{equation}
where $u(t)=(u_1(t),u_2(t)):[0,T]\rightarrow U.$

The cost functional in \eqref{cost} is: 
\begin{equation}\label{cost:1}
\ell(x,u,t) = \chi(x)_{B(0,10^{-4})},\qquad g(x) =\chi(x)_{B(0,10^{-4})}, \qquad \lambda = 0,
\end{equation}
where we measure the cost to reach a ball around the origin.
The corresponding HJB equation is a well-known time dependent eikonal equation. To show the quality of our approximation we compare our method with a classical approach based on a structured grid and state space interpolation (\cite{CF96}) using two different types of constraints. We will set $h = 0.005$ for all simulations in this test. The pruning is chosen according to the error estimates in \cite{SAF19} as $\ep=h^2$ to keep the first order of convergence of the value function. The discrete controls used to build the tree and to compute the value function are taken on the square $[-1,1]^2$ as:
$$\left\{(1, 0), (1, 1); (0, 1); (-1, 1); (-1, 0); (-1, -1); (0, -1); (1, -1); (0, 0)\right\}.$$ 
The choice of the square $[-1,1]^2$ for the controls will allow an efficient pruning, obtaining a tree which is equivalent to a grid constructed on the constrained domain. The initial condition will be $x = (1,1)$ and the final time $T=2$.
\paragraph{\bf Test 2a: the circular channel}
With the dynamics defined in \eqref{exact} we will set the domain $$\Omega=[-tol,1]^2 $$ and the following functions:

$$\varphi_1(x) =((x-1)^2+y^2-1.1;\quad \varphi_2(x) = -((x-1)^2+y^2-0.9);$$
are used to define the set $O_j, j=1,2$. The dynamics is then constrained to remain inside the circular channel $\Omega\setminus (O_1\cup O_2)$. To discretize this geometry we require a very fine spatial discretization with a classic method (see \cite{CF96}) whereas TSA only requires to follow the dynamics. In the classical method we set $\Delta x = 0.0025$ since the constraint is curvilinear and it needs a fine discretization to get a proper approximation. To build the optimal trajectories for both methods we consider the control $(0,0)$ and $64$ controls equidistributed on the ball centered in $(0,0)$ with radius $1$.

\begin{figure}[htbp]
\includegraphics[scale=0.4]{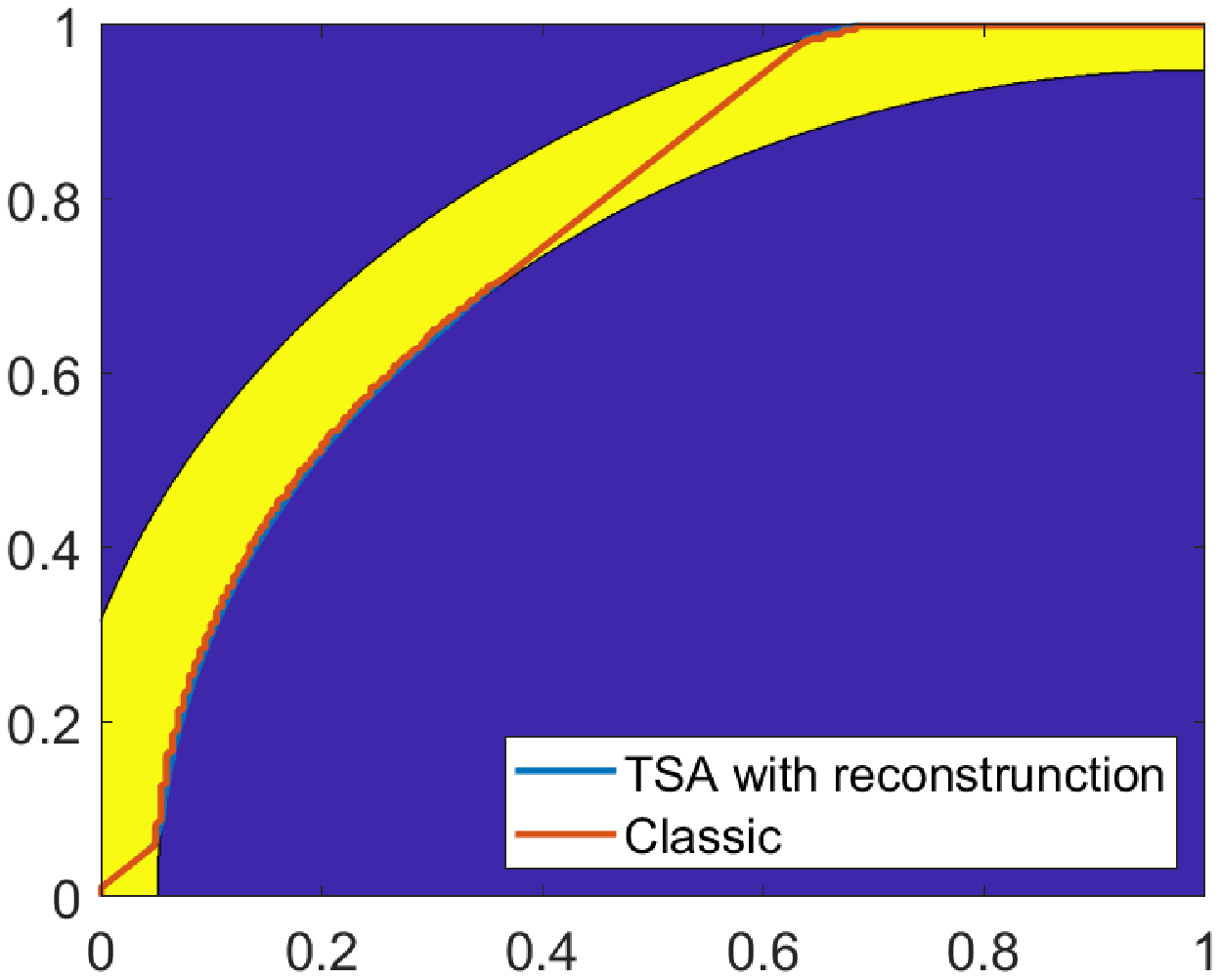}
\includegraphics[scale=0.41]{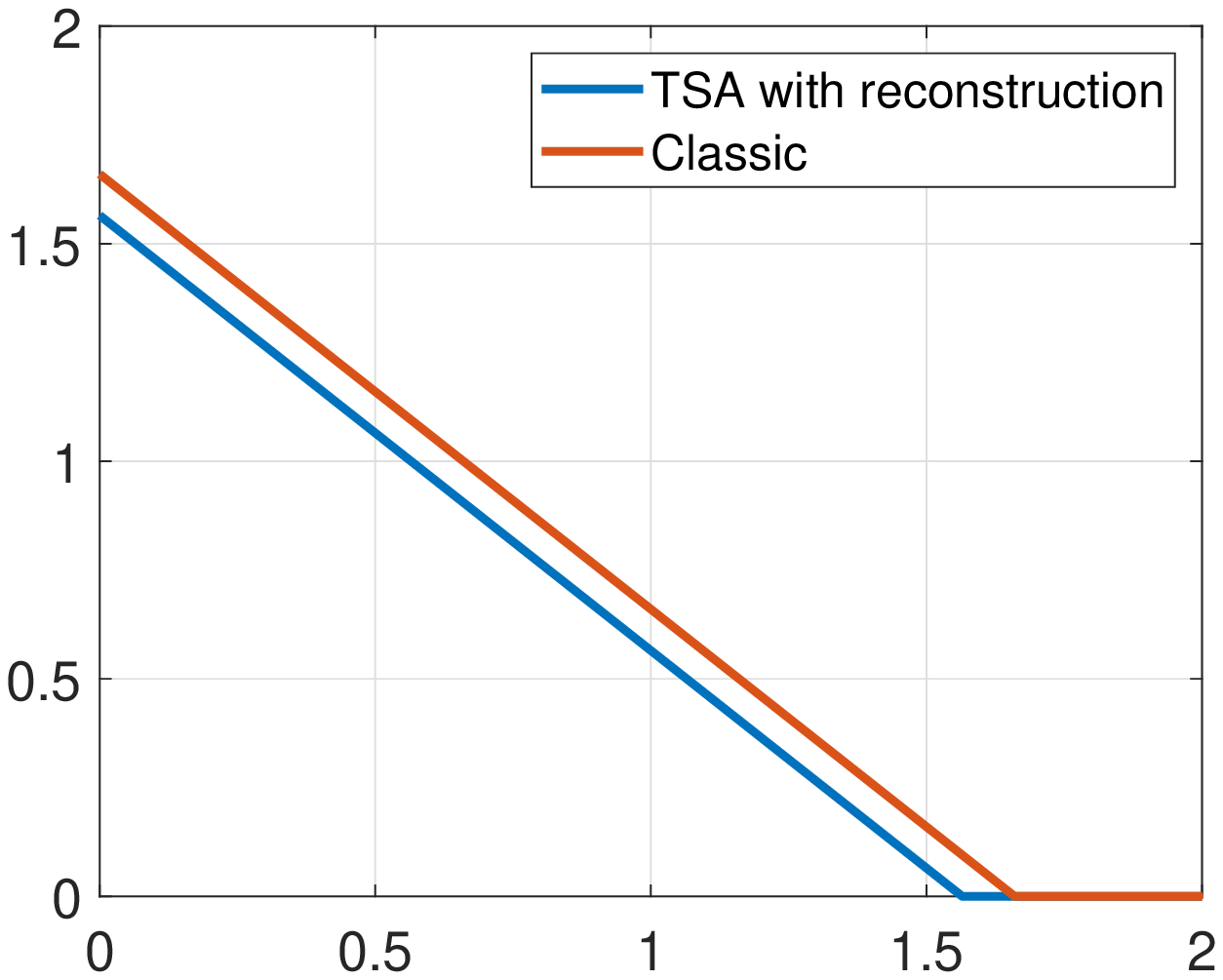}
\caption{Test 2a: Optimal trajectory (left) and cost functional (right).}
\label{fig1}
\end{figure}

In the left panel of Figure \ref{fig1}, we show the optimal trajectory using TSA method and the classical algorithm. We can easily see that the solutions are very similar as confirmed by the values of the cost function on the right panel of Figure \ref{fig1}. On the other hand we would like to mention that the TSA took about 7 seconds to compute the value function the feedback law whereas 14 seconds with the classical method.

\begin{figure}[htbp]
\includegraphics[scale=0.3]{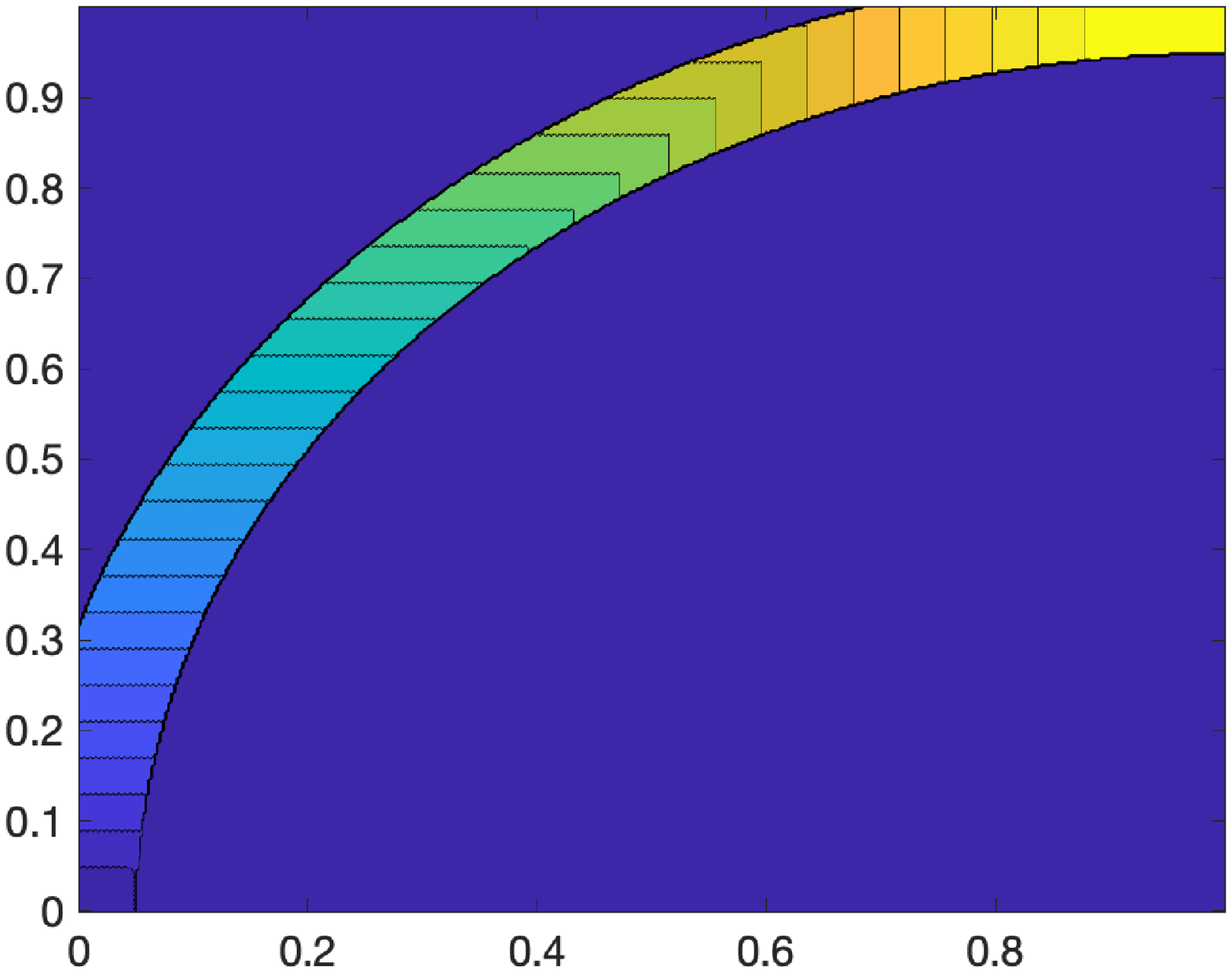}
\includegraphics[scale=0.3]{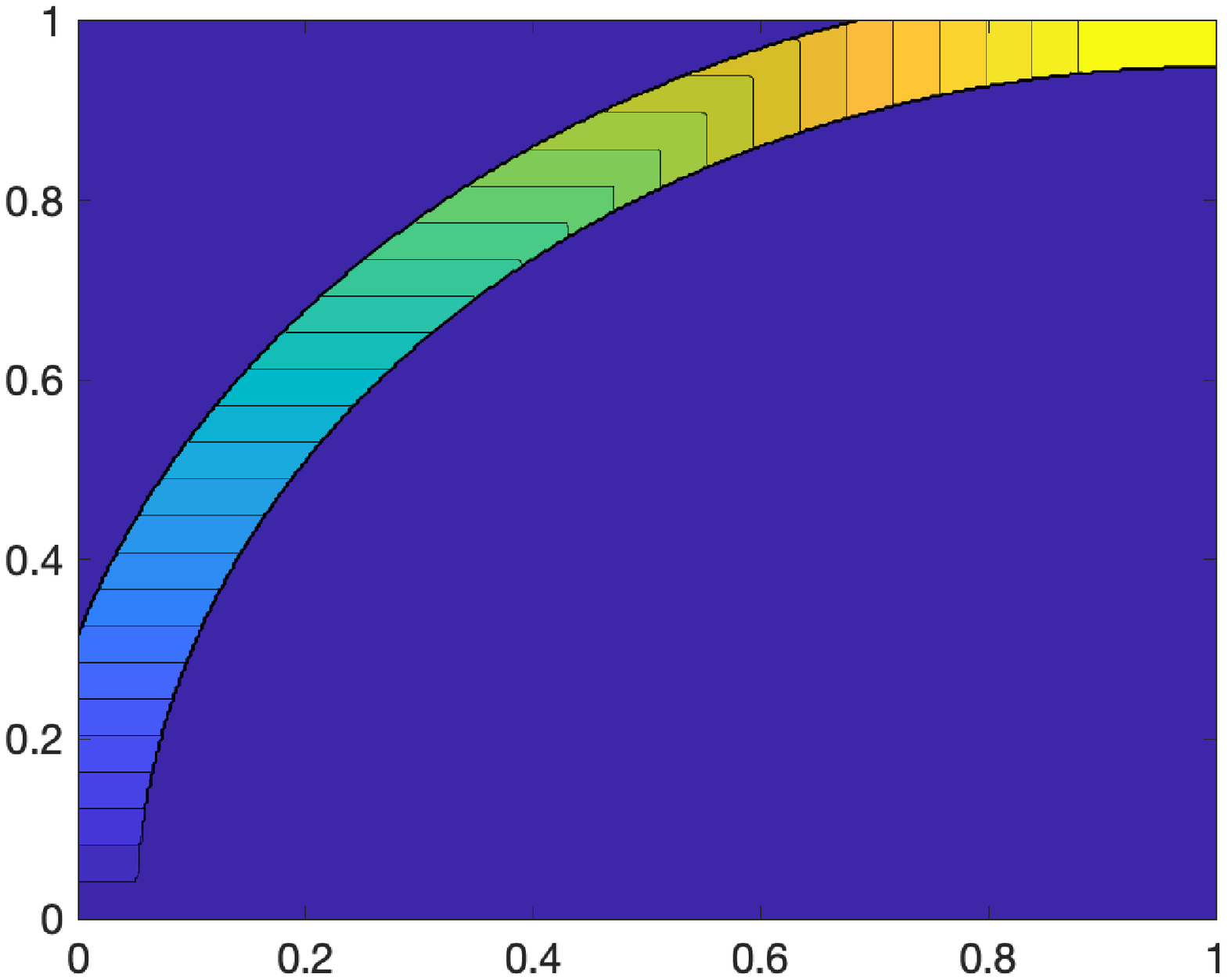}
\caption{Test 2a: Value function at time $t=0.75$ for the classical method (left) and for the TSA (right).}
\label{fig2}
\end{figure}

We show the contour lines of the value function at time $t=0.75$ with a classical method (left panel of Figure \ref{fig2}) and the TSA (right panel of Figure \ref{fig2}). It is clear that the pictures agree.

\begin{figure}[htbp]
\includegraphics[scale=0.4]{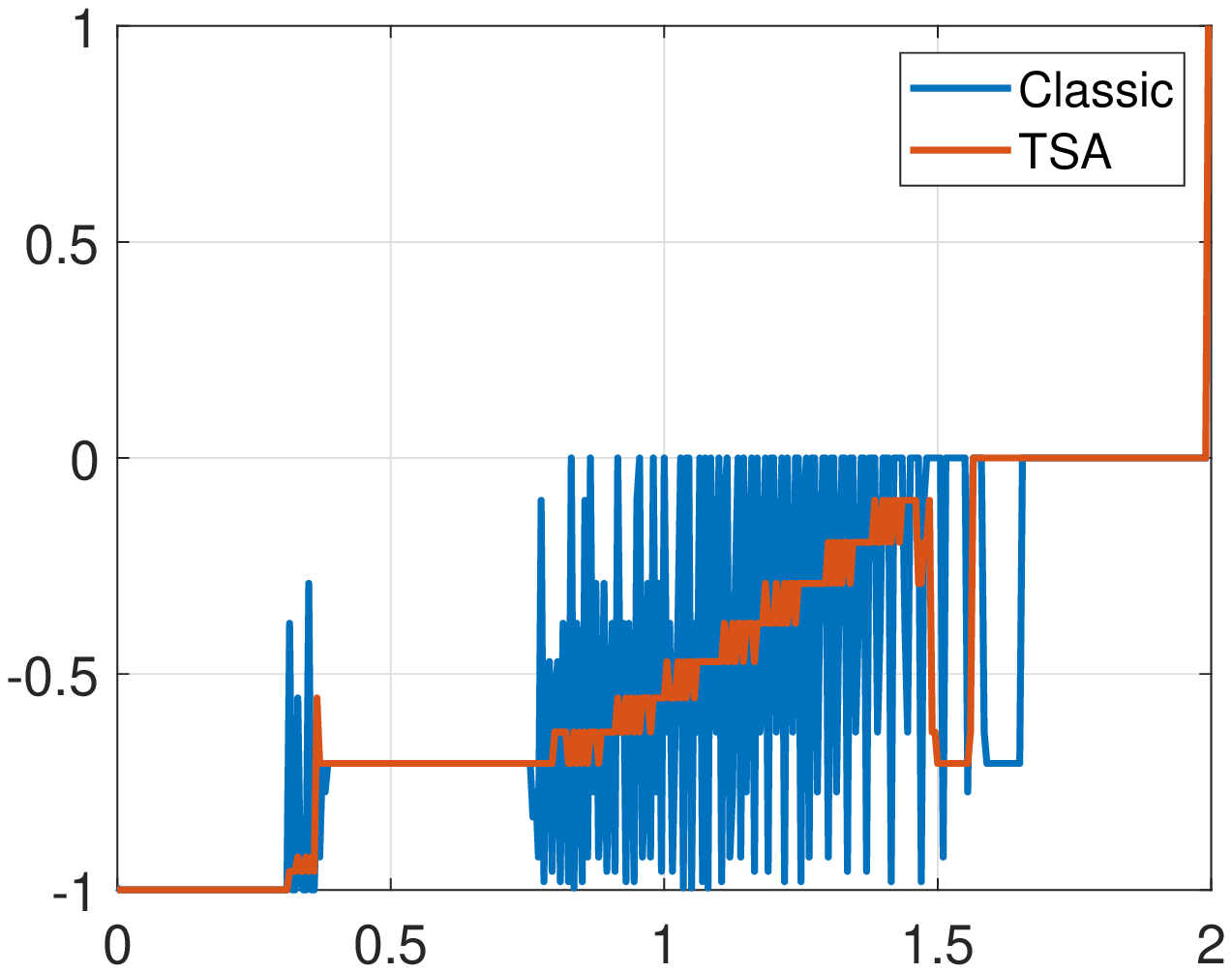}
\includegraphics[scale=0.4]{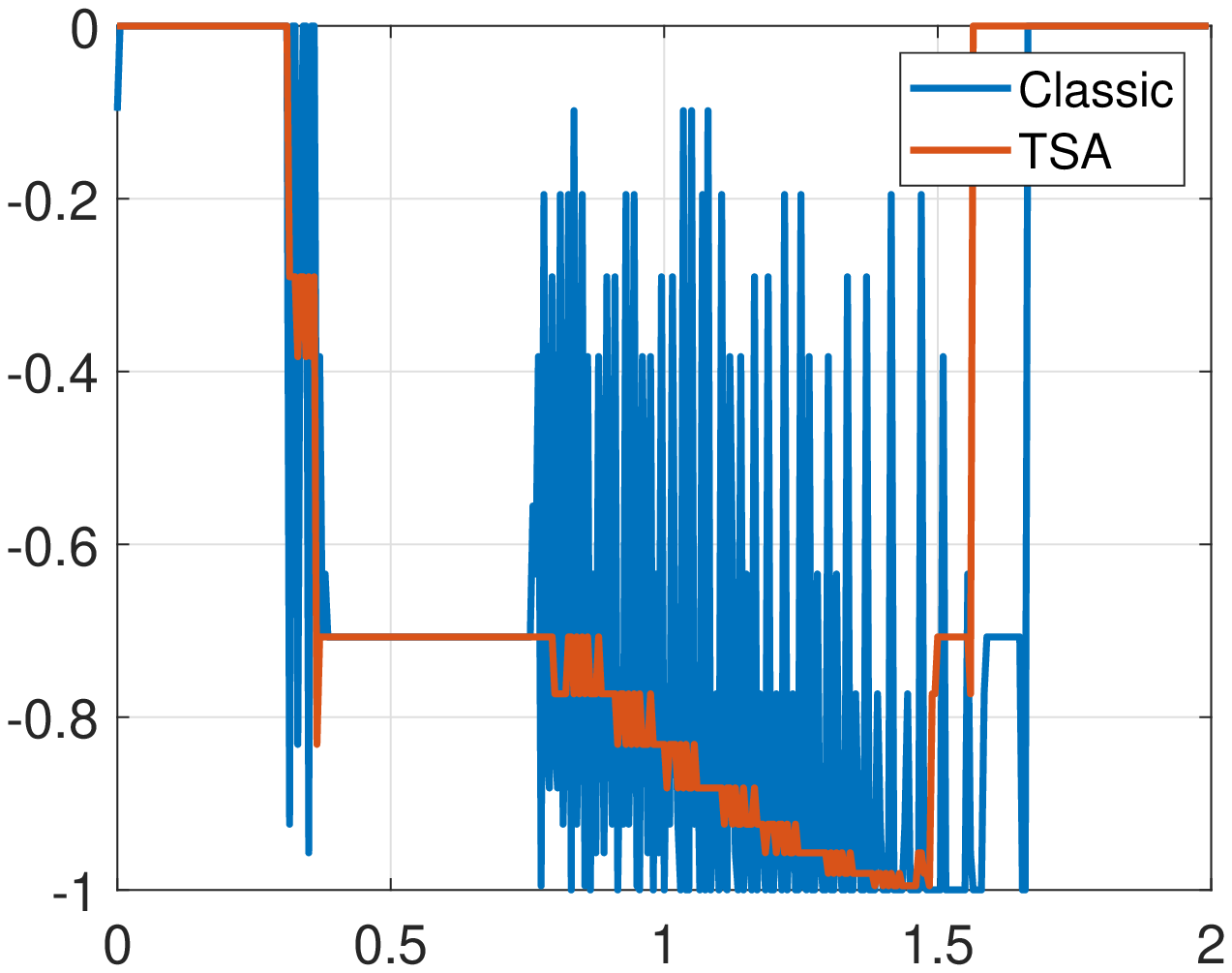}
\caption{Test 2a: First component (left) and second component (right) of the optimal control.}
\label{fig3}
\end{figure}

The optimal policy is shown in Figure \ref{fig3}. As one can see the controls have a high chattering behaviour which makes hard to reproduce this policy. This happens since our control space is not continuous and the control jumps between different values to reach the desired configuration. We also remind that, although the value function is unique, the optimal control is not. One can see that the control computed by TSA is different with respect to the control computed by the classical approach, but they lead to similar trajectories.

\begin{figure}[htbp]
\includegraphics[scale=0.4]{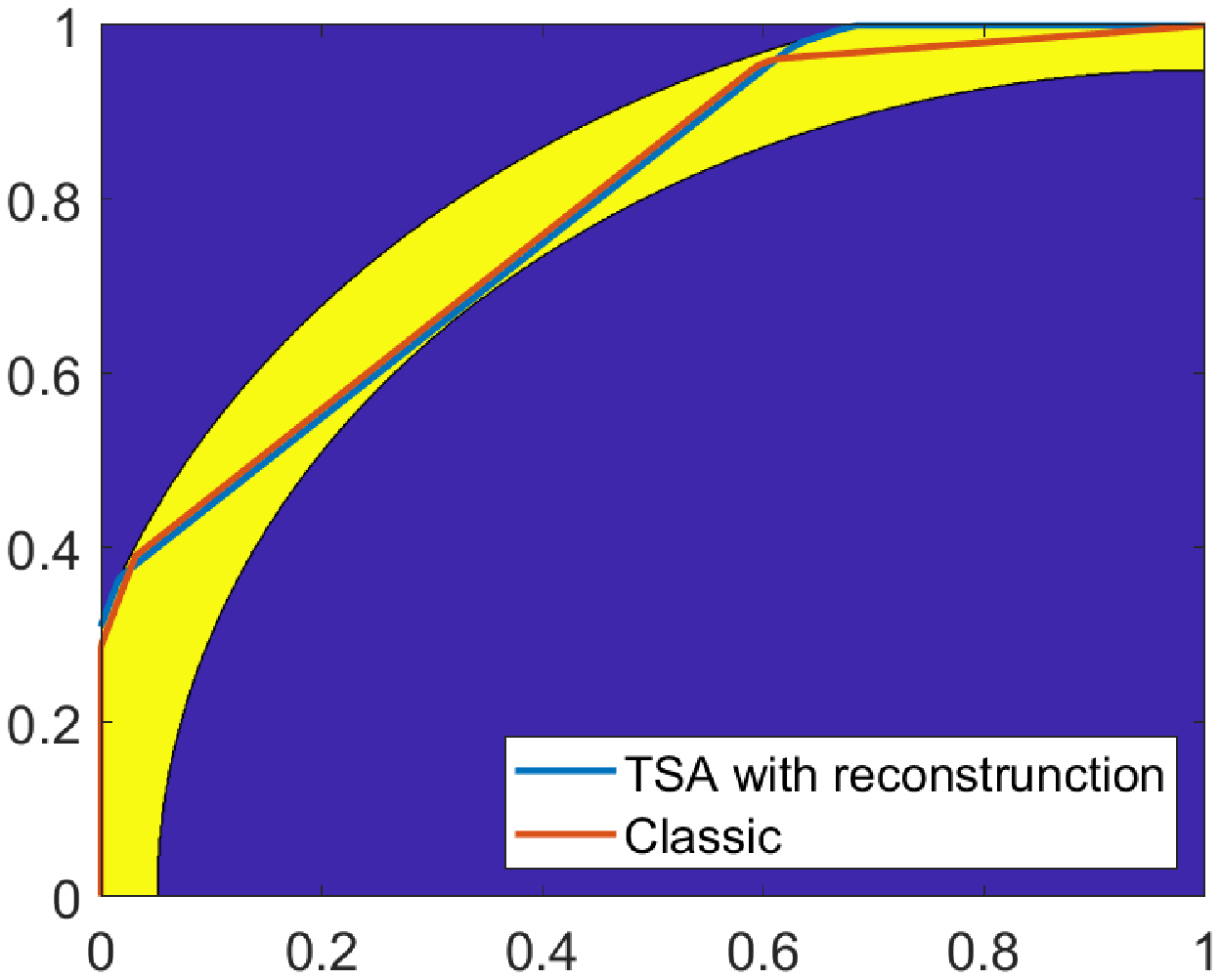}
\includegraphics[scale=0.4]{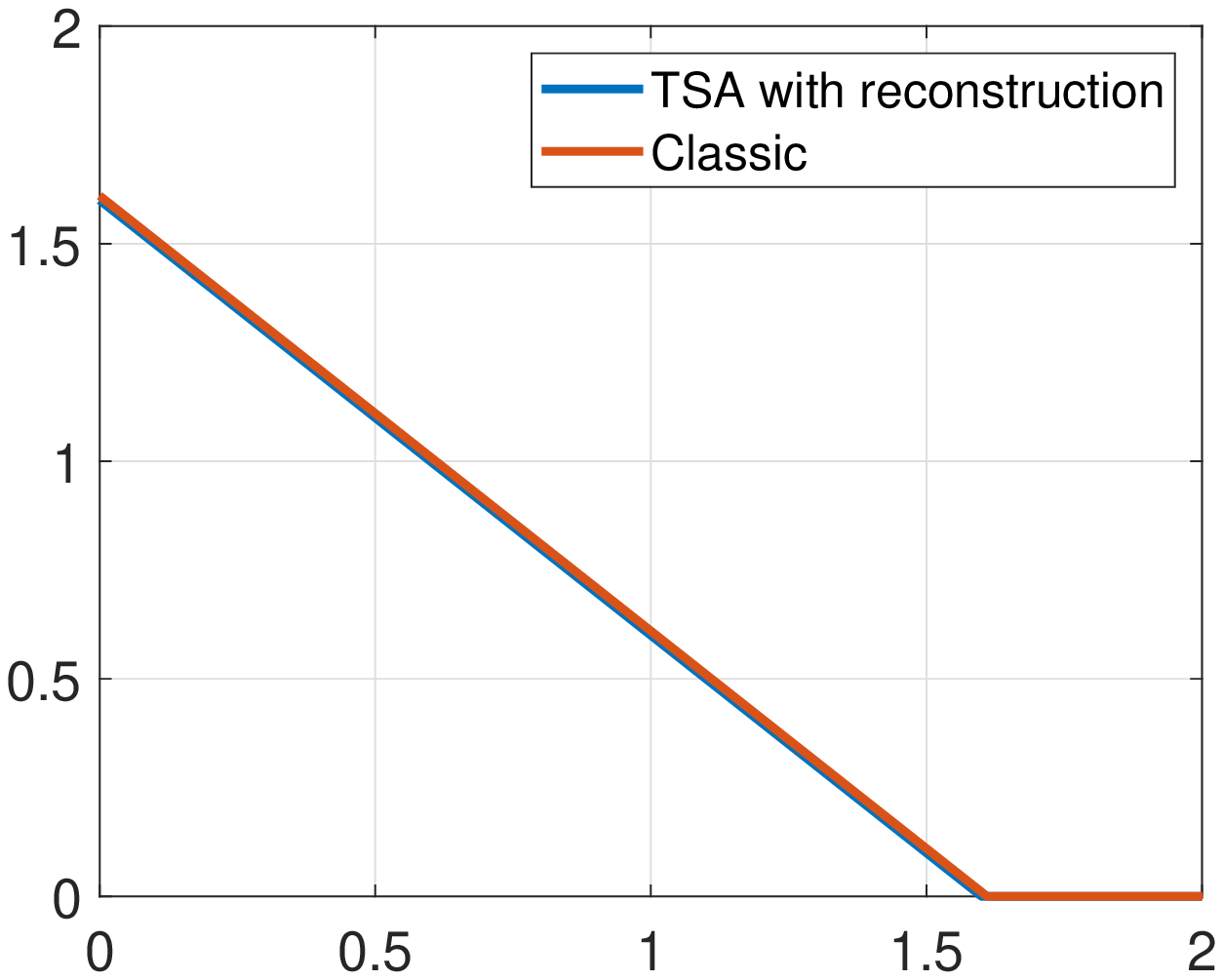}
\caption{Test 2a: Optimal trajectory (left) and cost functional (right) without chattering in the controls.}
\label{fig4}
\end{figure}

One can also introduce an inertia criteria in the feedback reconstruction to stabilize the feedback control as in \cite{F01} or \cite[Chap 8]{FF13}. This is often required for engineering applications and works penalizing at time $t_{i+1}$ a control which is far from the previous one computed at time $t_i$. We are going to consider the following feedback reconstruction for $n\ge 1$:
$$
u_{n}^{*}:=\argmin\limits_{u\in \widetilde{U}} \left\{V^{n+1}(\zeta^n_*+ h f(\zeta^n_*,u,t_n)) + h \,( \ell(\zeta^n_*,u,t_n) + \gamma |u-u^*_{n-1}|^2) \right\},
$$
where we added the term $ h \,\gamma |u-u^*_{n-1}|^2$ to reduce the distance between the new reconstruction and the previous control. We are going to fix $\gamma=7$.
 The optimal trajectory with this criteria is shown in the left panel of Figure \ref{fig4}. One can see that now the trajectories are a bit different from the previous ones.
\begin{table}[htbp]
		\centering
		\begin{tabular}{ccc}
					\hline
					  & With chattering & Without chattering  \\
					\hline
					TSA & 1.565 & 1595 \\
					Classic & 1.660 & 1.600 \\	
		\end{tabular}
   	\caption{Test 2a: Values of $V^0(x)$ for the two methods with different feedback reconstructions.}
	\label{test1a:table}
	\end{table}
We show in Table \ref{test1a:table} the value function at the initial time for the two methods with the two different feedback reconstructions. The classical approach reaches an improvement both in the stability of the feedback and in the cost functional, while the TSA gets a worse cost functional, since the previous reconstruction was not presenting a high chattering behaviour. For completeness we also show the nodes of the tree in the left panel of Figure \ref{fig4}. One can see that the nodes follow the constraint naturally without imposing any further restriction.

\begin{figure}[htbp]
\includegraphics[scale=0.4]{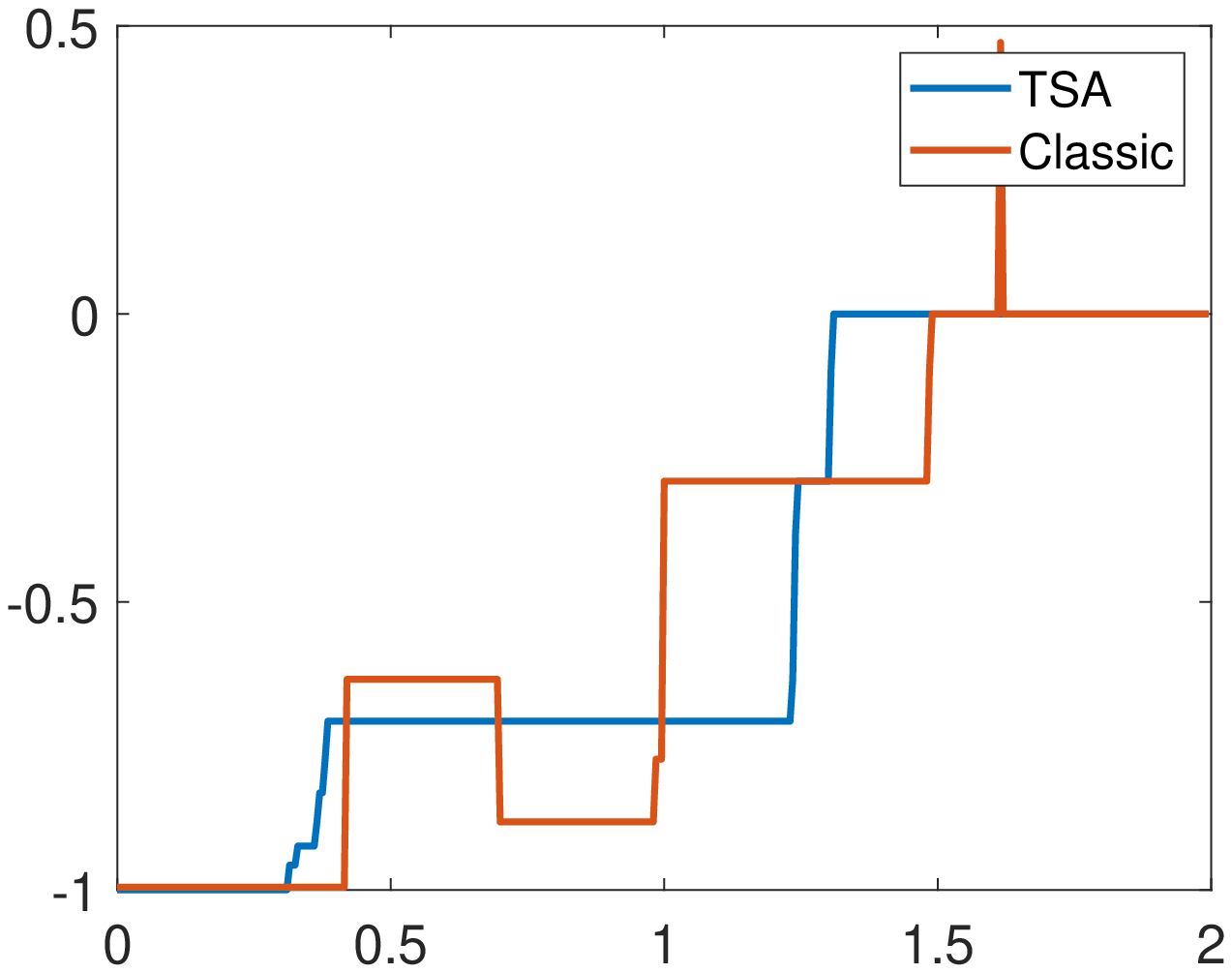}
\includegraphics[scale=0.4]{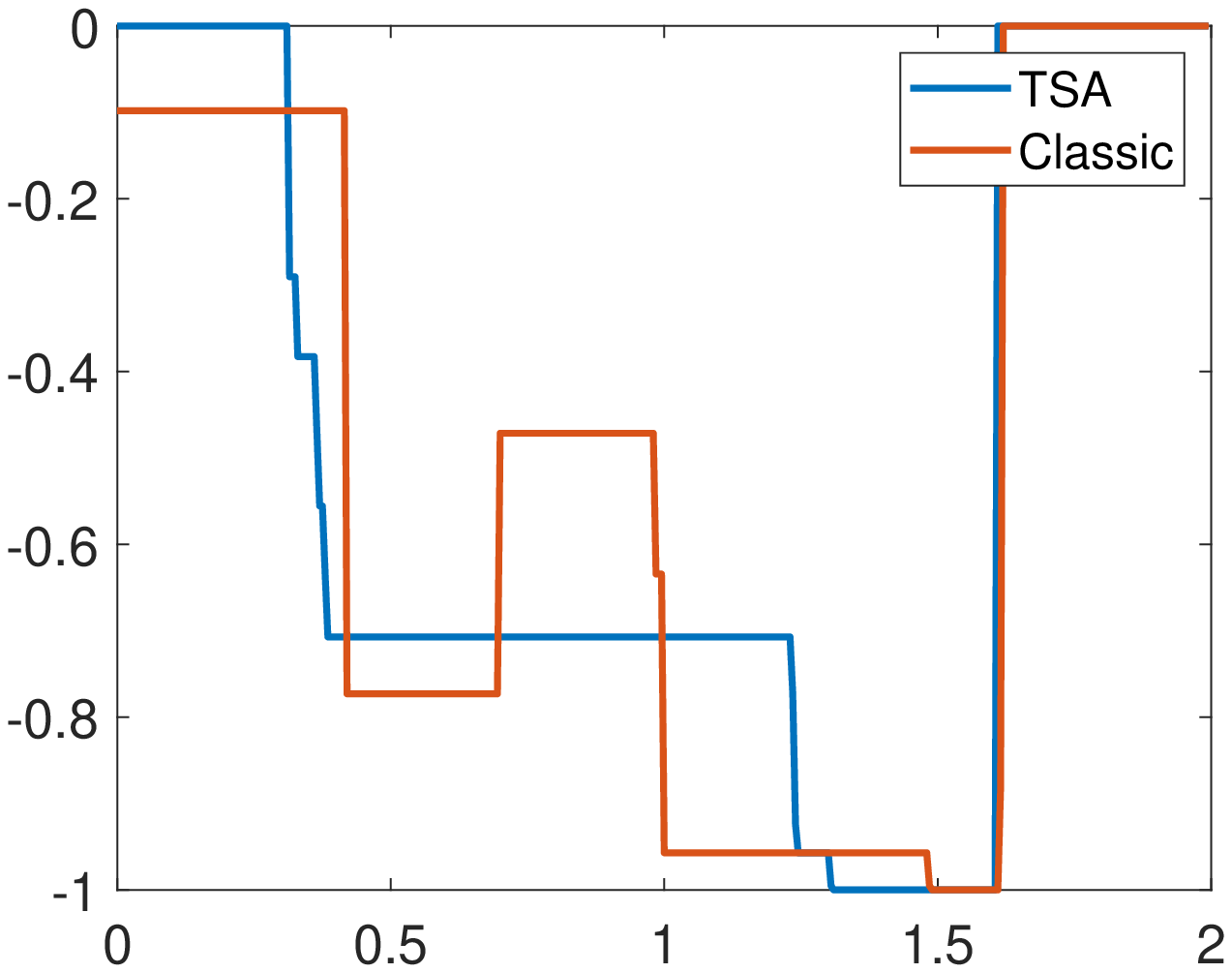}
\caption{Test 2a: First component (left) and second component of the optimal control without chattering.}
\label{fig5}
\end{figure}

The (sub)optimal controls without chattering are then shown in Figure \ref{fig5}. Now, the controls are rather stable and much more suitable for applications.

\paragraph{\bf Test 2b: a channel with obstacles}
We modify our constraint using, again, the dynamics defined in \eqref{exact} and the running cost \eqref{cost:1}. We will set the domain as the yellow part in the left panel of Figure \ref{fig3} with the obstacles 

$$\varphi_1(x,y) = -((x-0.9)^2+(y-0.9)^2)-0.005)$$
$$ \varphi_2(x,y) =  -\left(\dfrac{(x-0.3)^2}{80}+(y-0.05)^2-0.001\right)  $$

%

\begin{figure}[htbp]
\includegraphics[scale=0.4]{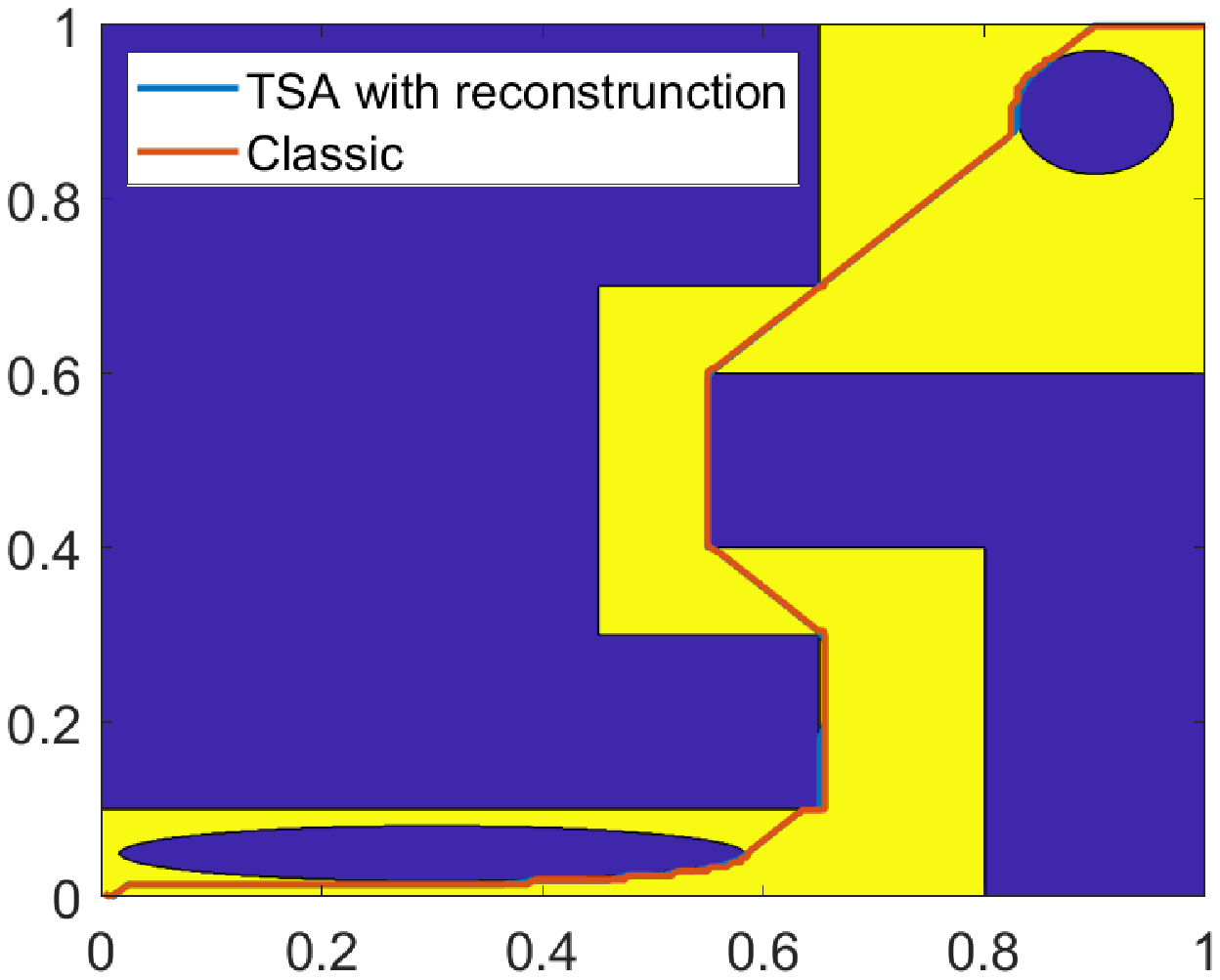}
\includegraphics[scale=0.42]{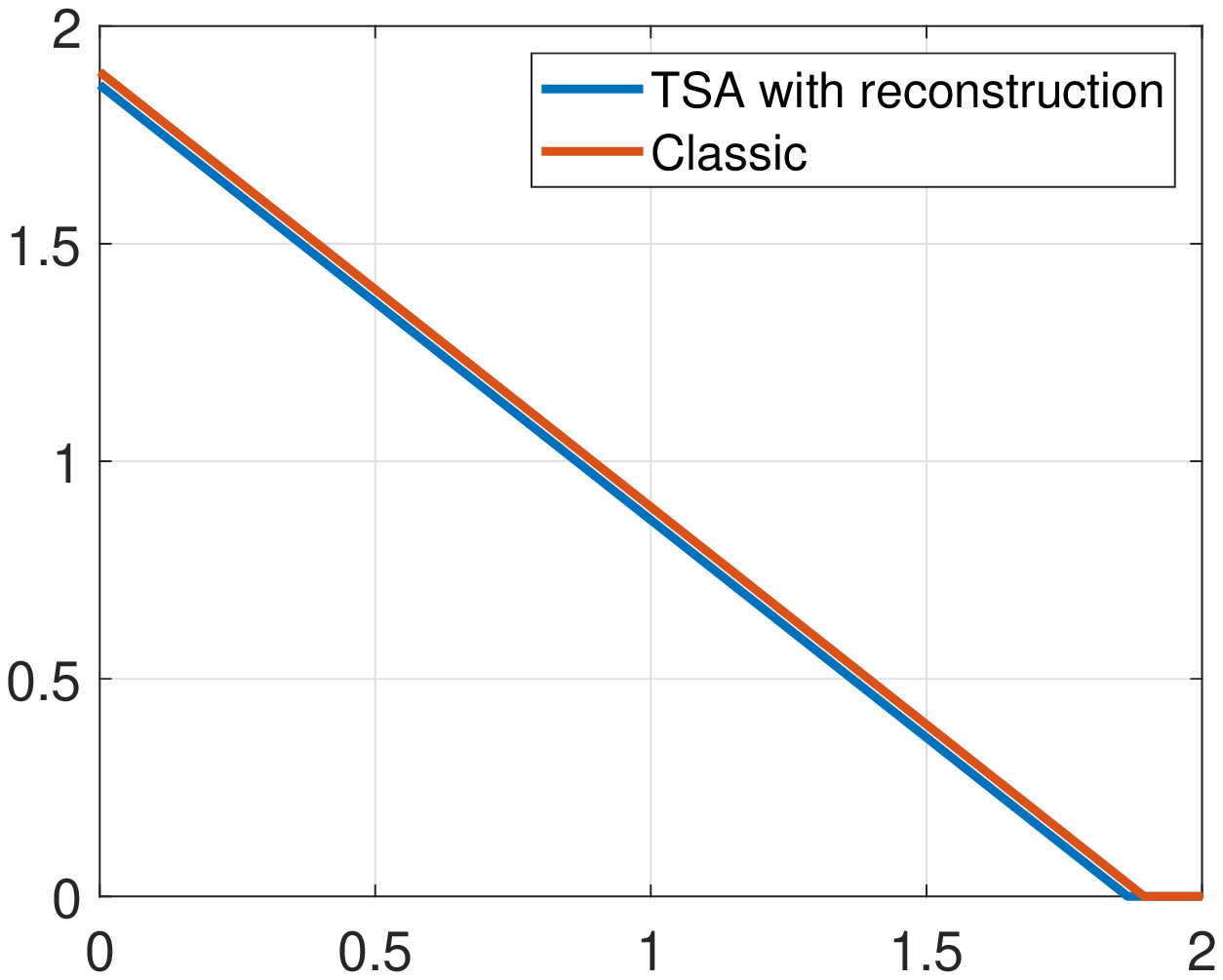}
\caption{Test 2b: Optimal trajectory (left) and cost functional (right).}
\label{fig6}
\end{figure}

The left pannel of Figure \ref{fig6} shows the optimal trajectory. We can see that it is hard to distinguish between the solution driven by the classical approach and the TSA method. They both try to avoid the first obstacle to reach fastly the closest corner and to continue along the border till the origin avoiding also the elliptical constraint. The evaluation of the cost functionals is very similar (see the right panel of Figure \ref{fig6}). However, the CPU time is 17s with the TSA and 25s with a classical approach, this means a reduction of about the 33\%. In this setting we use the same parameters of Test 1a except for the number of controls in the reconstruction which are $32$ plus the origin.

\begin{figure}[htbp]
\includegraphics[scale=0.4]{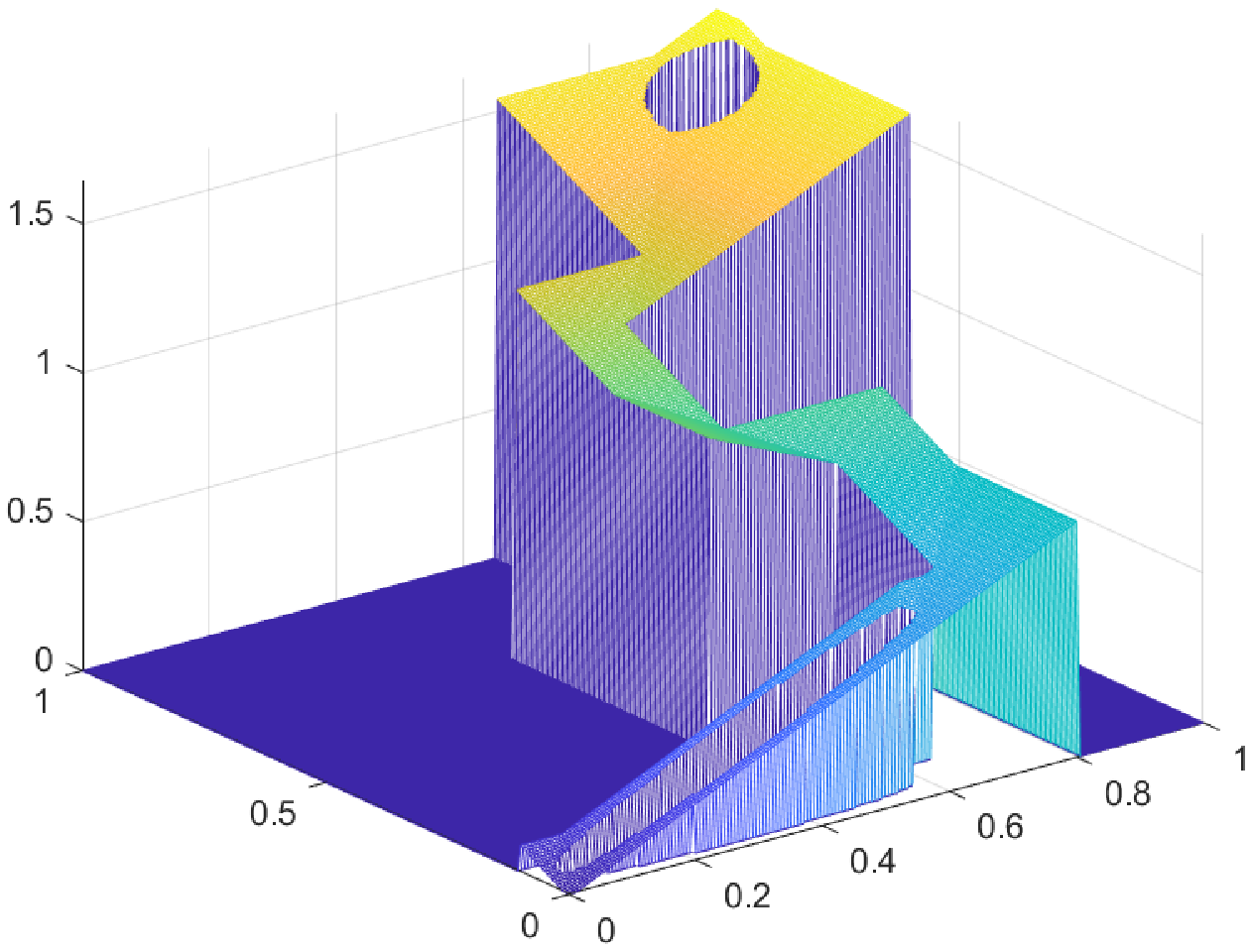}
\includegraphics[scale=0.4]{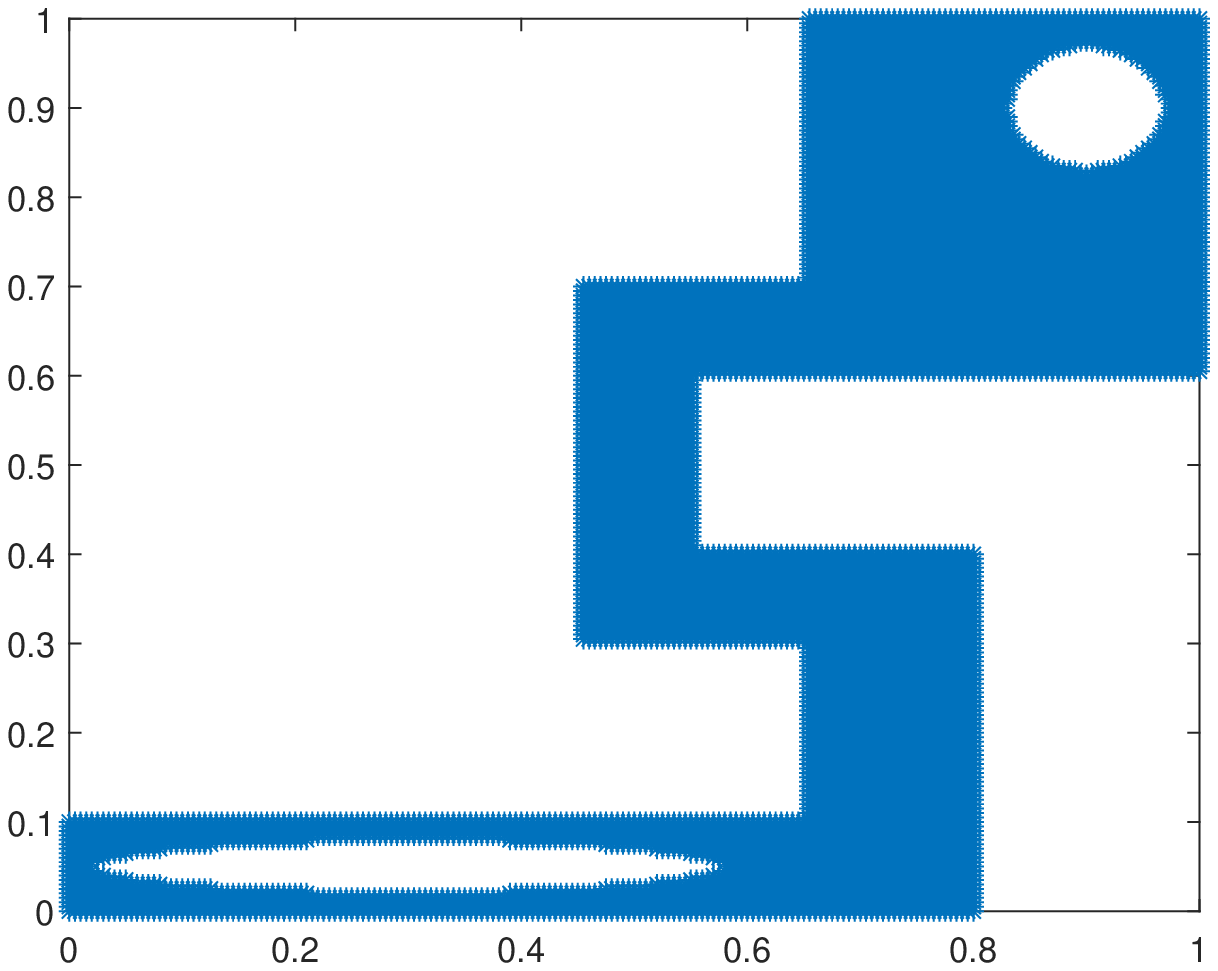}
\caption{Test 2b: Value function at time $t=0$ (left) and nodes of tree (right).}
\label{fig7}
\end{figure}

The value function at the initial time is plotted in the left panel of Figure \ref{fig7}. We also show the nodes of the tree in the right panel. We again want to emphasize how it is easy with a TSA method to remain inside the constraint due to the pruning criteria which also involves the state constraints $\zeta_i^n\in O_j$.

\begin{figure}[htbp]
\includegraphics[scale=0.4]{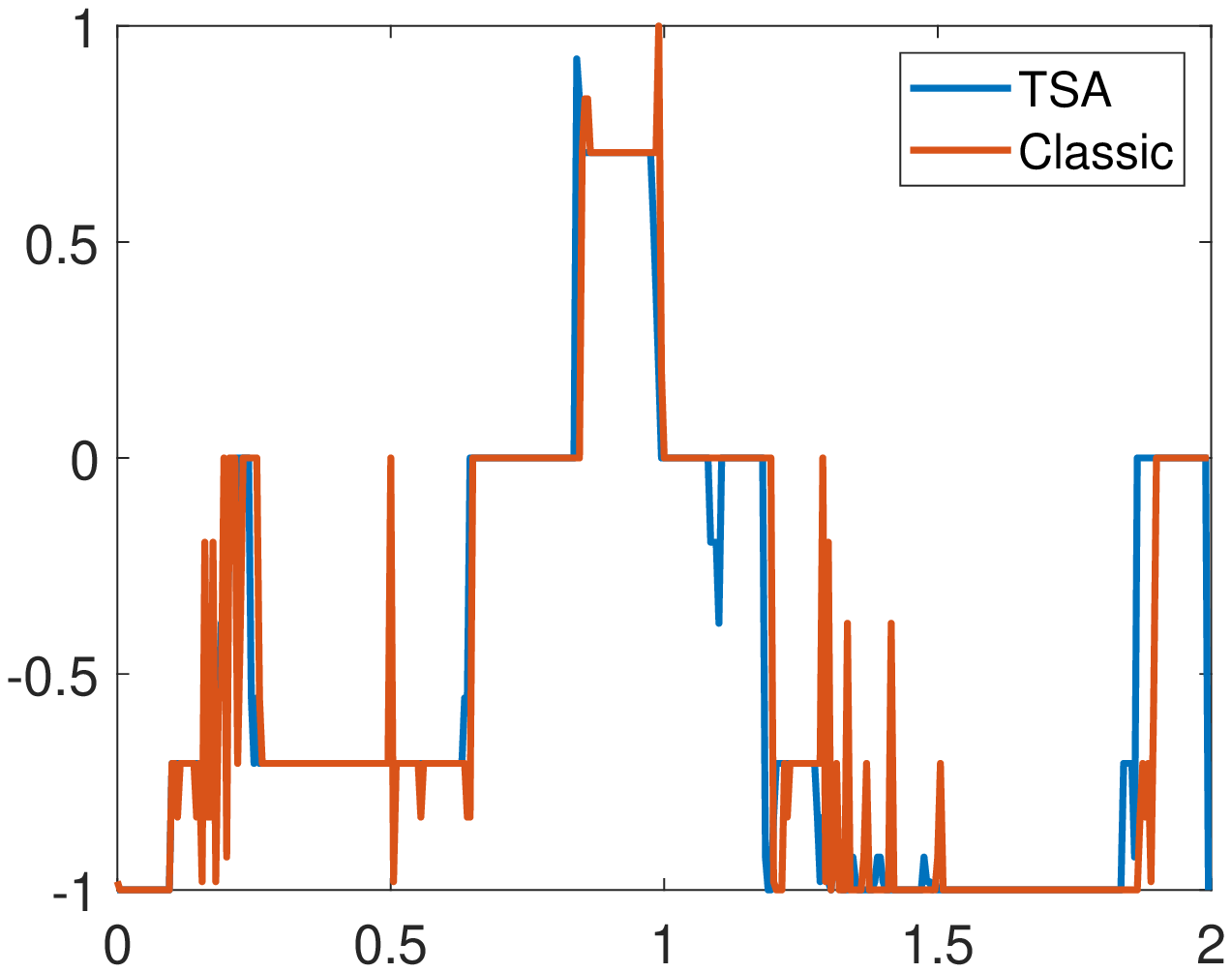}
\includegraphics[scale=0.4]{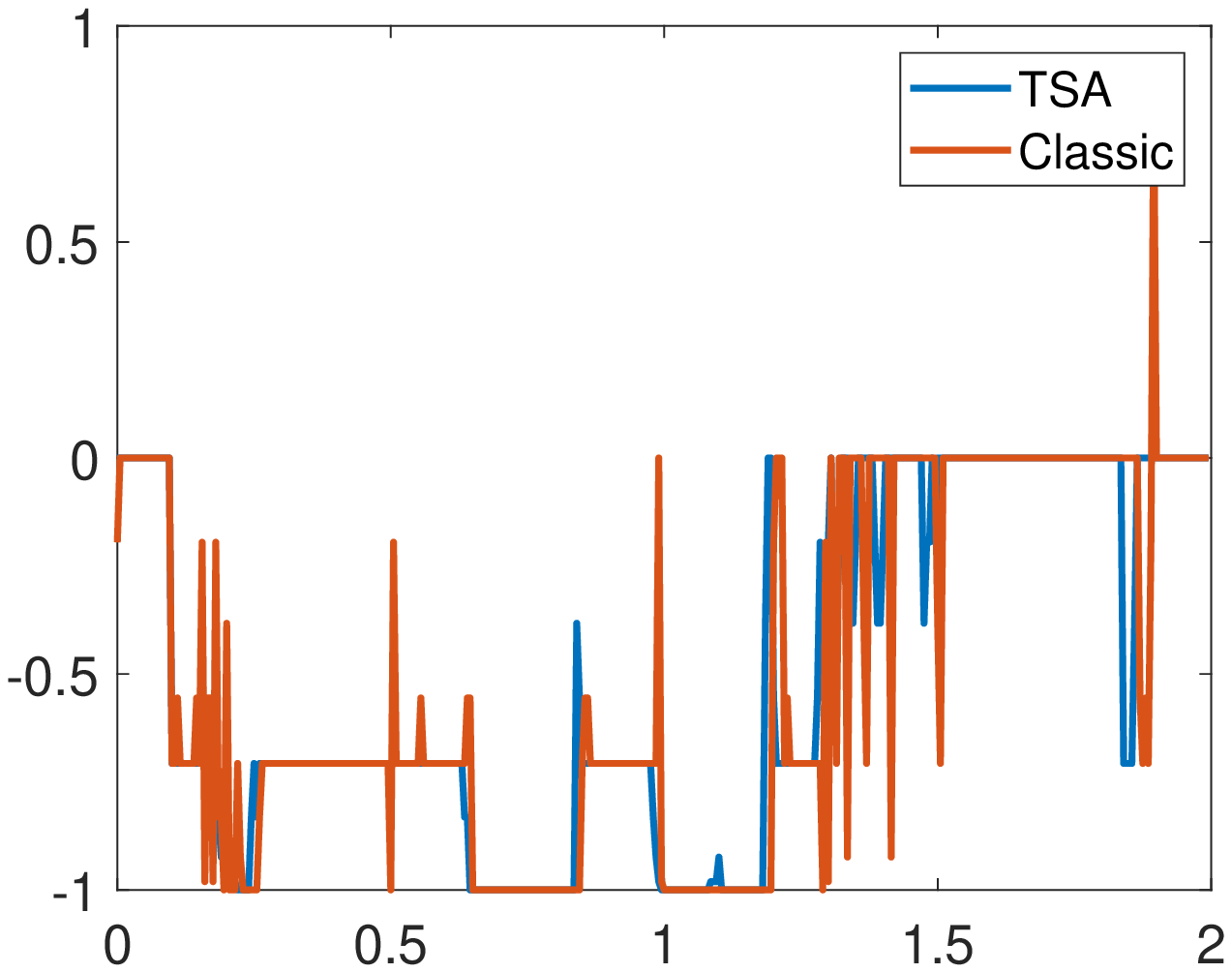}
\caption{Test 2b: First component of the optimal control (left) and second component of the optimal control.}
\label{fig8}
\end{figure}

We finally show the optimal policy for this problem in Figure \ref{fig8}. We can see that policies show a good agreements between the methods.

\subsection{Test 3: Constrained Van der Pol}
In this last test we consider the Van der Pol oscillator. The dynamics in \eqref{eq} is given by
\begin{equation}\label{eq:vdp}
	f(x,u)=
		  \begin{pmatrix}
         x_2 \\
         0.15(1-x_1^2)x_2-x_1+u
          \end{pmatrix}\quad u\in U\equiv [-1, 1].
\end{equation}
The cost functional in \eqref{cost} is: 
\begin{equation}
\ell(x,u,t) = \|x\|^2 \qquad g(x) = \|x\|^2, \qquad \lambda = 0,
\end{equation}
and aims to steer the solution to the origin which is a repulsive point for the uncontrolled dynamics $(u=0)$. We set $T=1.4$, $x=(0.4,-0.3)$ as initial condition, $h=0.025$ and $\ep = h^2$. Furthermore, the constraint is the box: $$\Omega:=(-h, 0.5) \times (-\infty, 0.1) \setminus \{(0.1,0.3)\times (-0.5,-0.3)\}$$

\begin{figure}[htbp]
\includegraphics[scale=0.3]{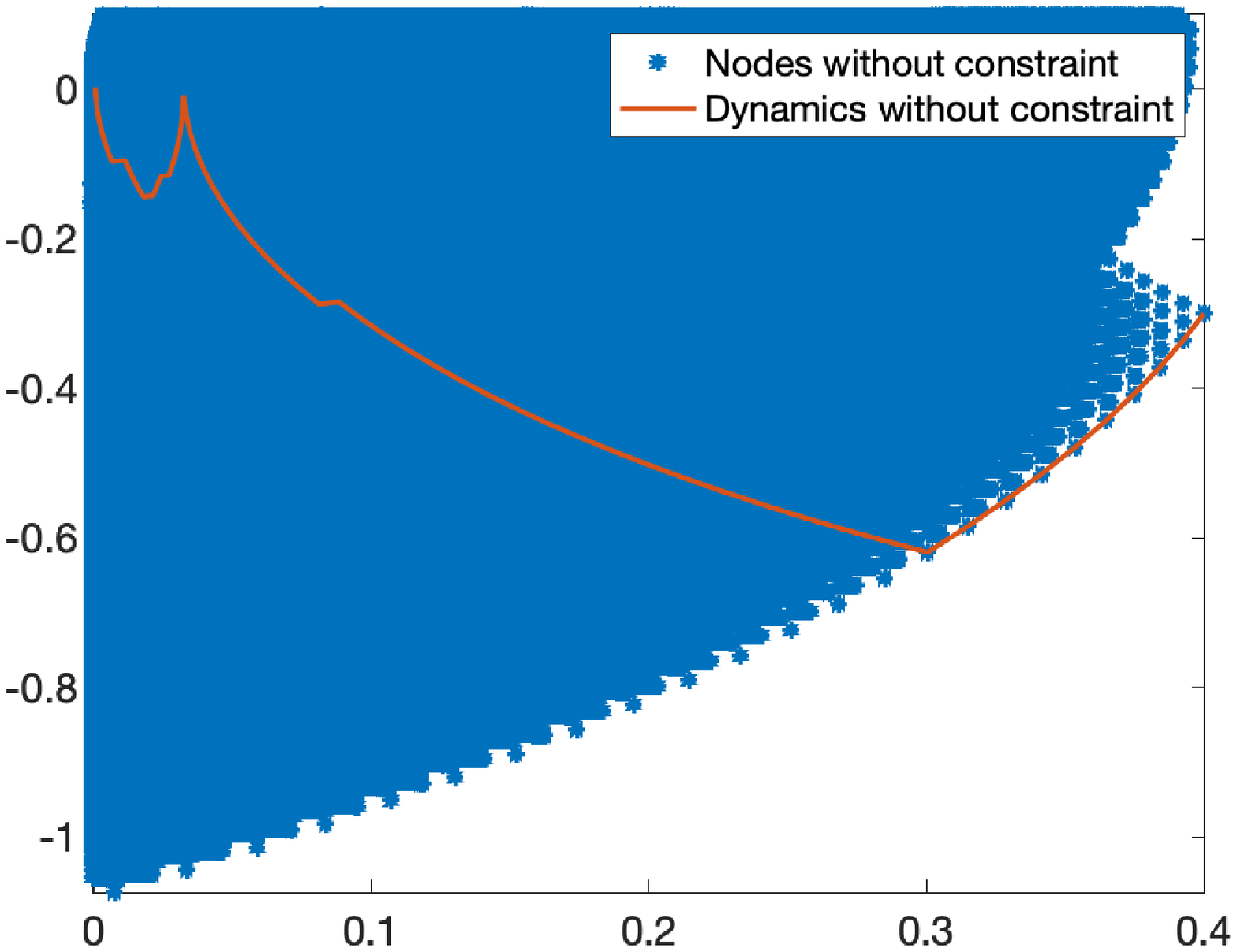}
\includegraphics[scale=0.3]{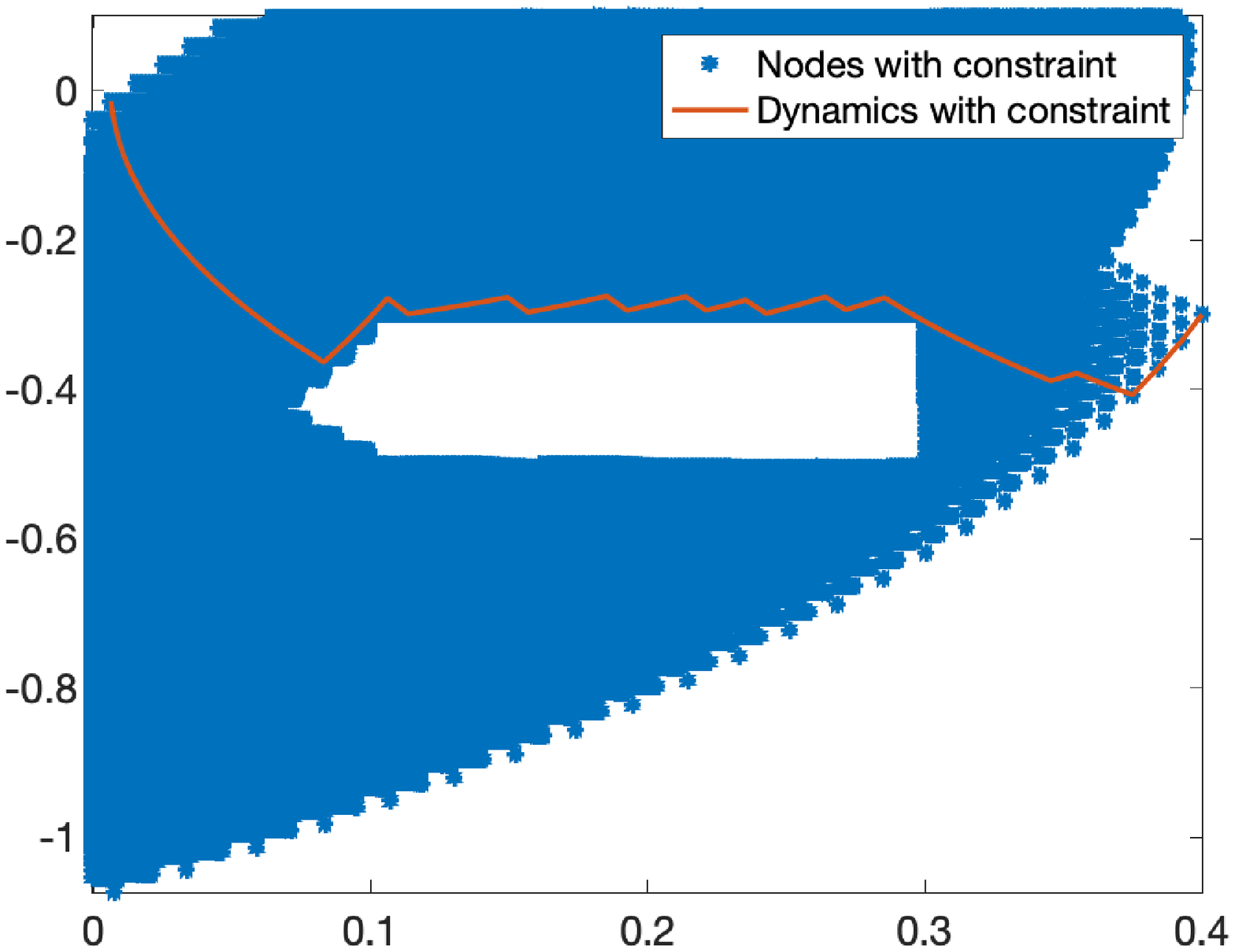}
\includegraphics[scale=0.35]{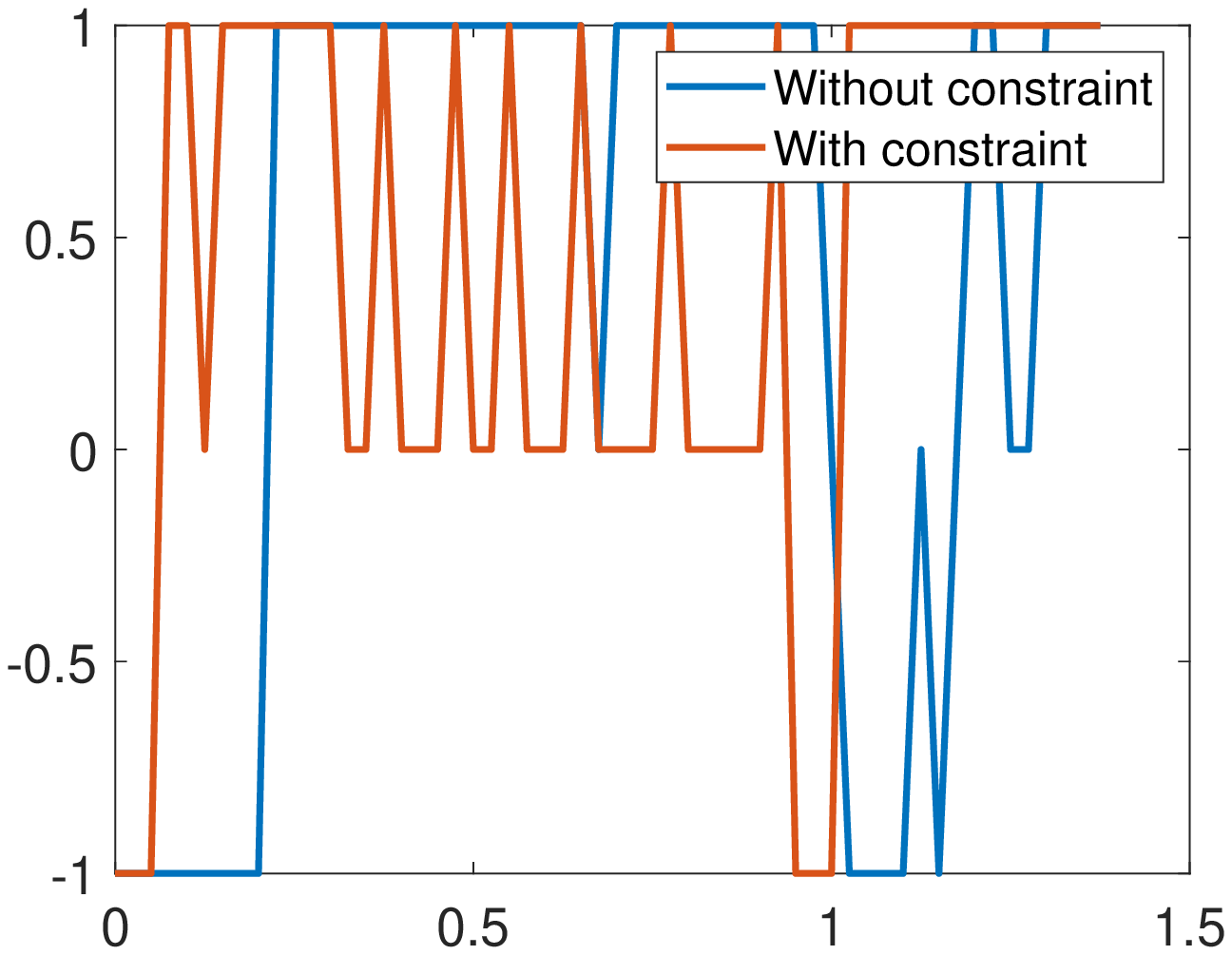}
\caption{Test3: Top: Optimal trajectory with constraints (left) and  optimal trajectory without constraints (right). Bottom: optimal control.}
\label{fig9}
\end{figure}

To compute the value function we use two discrete controls $\{-1,1\}$, whereas for the feedback control we use $3$ discrete controls: $\{-1,0,1\}$.
The optimal trajectories are shown in the top panel of Figure \ref{fig9}. In this example we compare the optimal trajectory with and without constraints. We can see in the left panel of Figure \ref{fig9} that the trajectory is passing through the constraint and also the tree nodes covers that region, whereas in the right panel the trajectory avoids the obstacle.
In this figure we can also see not only the constraint (the rectangle), but also other spacial regions not reachable by the constrained dynamics.
 In this way we can see how the trajectories behave differently when a constraint is added. At the end they both reach the origin. To obtain those trajectories we have computed two different policies as shown in the bottom panel of Figure \ref{fig9}.


\section{Conclusions}
We have examined some optimal control problems with state constraints from the numerical point of view. In the first part we have  given a new formulation of the infinite horizon problem  with convex constraints and we have proved a convergence result for a discrete time approximation  using also tools of multivalued analysis.

In the second part,  we have worked on finite horizon control problems proposing  an extension of the TSA in order to reduce the complexity of the DP algorithm. 
This approach does not need a fixed grid in space and exploits a tree structure for the approximation. We have also introduced   the synthesis of feedback controls by means a scattered interpolation routine which allows to increase the set of discrete controls in the reconstruction of optimal trajectories. In fact,  it is usually very common to compute the value function with a low number of controls and the feedback with a larger set of admissible controls to improve the accuracy. We have shown the effectiveness of the TSA by several numerical examples with convex and non-convex constraints (that are at present outside the limit of our convergence result). To this end we have compared our approach with a classical grid approach based on the interpolation on a fixed space grid. As for the case of optimal control problems without state constraints, the advantage of  the proposed TSA approach is its capability to deal with high dimensional problems, e.g. control of PDEs as shown in \cite{AFS18, AFS18b,  AS19}. We will address this problem in a future work.

%
%

\renewcommand\refname{References}

\end{document}